\documentclass[12pt,a4paper,twoside,final,notitlepage, leqno]{article}
\usepackage[english]{babel}
\usepackage[T1]{fontenc}  
\usepackage{epsfig, graphicx, amssymb}
\usepackage{amsmath,amsthm,epsfig,amsfonts,bbm}  
\usepackage{float}
\usepackage{color}
\def \smb {{\scriptstyle \bullet }}
\newcommand{\monitem}{ \smallskip \noindent $\bullet$ \quad  } 
\newcommand{\moneq}{\vspace*{-7pt} \begin{equation} \displaystyle } 
\newcommand{\moneqstar}{\vspace*{-6pt} \begin{equation*} \displaystyle } 
\newcommand{\monendstar}{\vspace*{-6pt} \end{equation*}   }
\newcommand{\monend}{\vspace*{-7pt} \end{equation}   }

\def\section*#1{}

\usepackage{fancyhdr}
\renewcommand{\headrulewidth}{0pt}

\begin{document} 

\fancypagestyle{plain}{ \fancyfoot{} \renewcommand{\footrulewidth}{0pt}}
\fancypagestyle{plain}{ \fancyhead{} \renewcommand{\headrulewidth}{0pt}}

~

\bigskip \bigskip 

\centerline {\bf \LARGE   Generalized bounce back boundary condition }

 \bigskip 
\centerline {\bf \LARGE  for the nine velocities  two-dimensional  }

\bigskip 
\centerline {\bf \LARGE  lattice Boltzmann scheme   }

 \bigskip \bigskip \bigskip \bigskip

\centerline { \large    Fran\c{c}ois Dubois$^{ab}$, Pierre Lallemand$^{c}$ and Mohamed Mahdi Tekitek$^{d}$}

\smallskip

\centerline { \it  \small   
$^a$   Dept. of Mathematics, University Paris-Sud,  B\^at. 425, F-91405  Orsay, France.} 

\centerline { \it  \small   
$^b$    Conservatoire National des Arts et M\'etiers, LMSSC laboratory,  F-75003 Paris, France.} 

\centerline { \it  \small  $^c$   Beijing Computational Science Research Center, 
Haidian District, Beijing 100093,  China.}

\centerline { \it  \small  $^d$   Dpt. Mathematics, Faculty of Sciences of Tunis, 
University Tunis El Manar, Tunis, Tunisia. } 

\bigskip      

\centerline {  05 July 2017 {\footnote 
    {\rm  \small $\,$ Contribution published in {\it Computers and Fluids},
volume 193, 103534, october 2019, \\  doi:  doi.org/10.1016/j.compfluid.2017.07.001.}} }

 \bigskip \bigskip 

 {\bf Keywords}: Taylor expansion method, acoustic scaling. 

  {\bf PACS numbers}:  
 02.70.Ns, 
  05.20.Dd, 
 47.11.+j. 
 
  {\bf AMS classification}: 
 76M28 

\bigskip \bigskip 
\bigskip \bigskip 

{\bf Abstract }

In a previous work \cite{DLT15}, we have proposed   
a method for the analysis of the bounce back boundary condition 
with the Taylor expansion method in the linear case. 
In this work two new schemes of modified bounce back are proposed. 
The first one is based on the expansion of the iteration 
of the internal scheme of the lattice Boltzmann method. 
The analysis puts in evidence some defects and a generalized version 
is proposed with a set of essentially four possible parameters to adjust. 
We propose to reduce this number to two with the elimination of 
 spurious density first order terms. 
Thus  a new scheme for bounce back is found exact up to second order and
allows an accurate simulation of the Poiseuille  flow for a specific
combination of the relaxation  and boundary coefficients.
We have validated  the general expansion of the value in the first cell 
in terms of given values on the boundary 
 for a stationary  ``accordion'' test case.


\newpage \bigskip {\bf Introduction}

In this contribution, we study boundary conditions for lattice Boltzmann schemes 
using the Taylor expansion
method proposed in our previous work \cite {DLT15}.
In that work, we have proposed a method for the analysis of the bounce back boundary condition
in the particular case of the D2Q9 scheme \cite {LL00} for a bottom boundary.
Note that the bounce back boundary condition and anti-bounce back boundary condition 
were studied by Ginzburg and Adler \cite{GA94} , Zou and He \cite{ZH97}, 
Bouzidi {\it{et al}} \cite{BFL01} and d'Humi\`eres and Ginzburg \cite{dHG03}. 
A particular choice of the LB parameters can enhance the precision of the scheme.

In this contribution, more general bounce back boundary conditions are proposed. 
We follow the same method as in \cite{DLT15}
to analyze the proposed scheme up to second order in space. 
Three schemes are investigated and  implemented  for a Poiseuille flow with
an imposed pressure field at the input and at the output of the domain.
Finally,  we propose a new scheme for bounce back exact up to order two in space
that allows an accurate simulation of the Poiseuille flow for any combination of the relaxation
coefficients.

\bigskip {\bf 1) \quad D2Q9 lattice Boltzmann scheme} 

\fancyhead[EC]{\sc{ F. Dubois,  P. Lallemand,  M.M. Tekitek }} 
\fancyhead[OC]{\sc{Generalized bounce back boundary condition for D2Q9 lattice Boltzmann}} 
\fancyfoot[C]{\oldstylenums{\thepage}}

\begin{figure}[htbp!]         
\begin{center}
\includegraphics[width=.40 \textwidth, angle=0]{{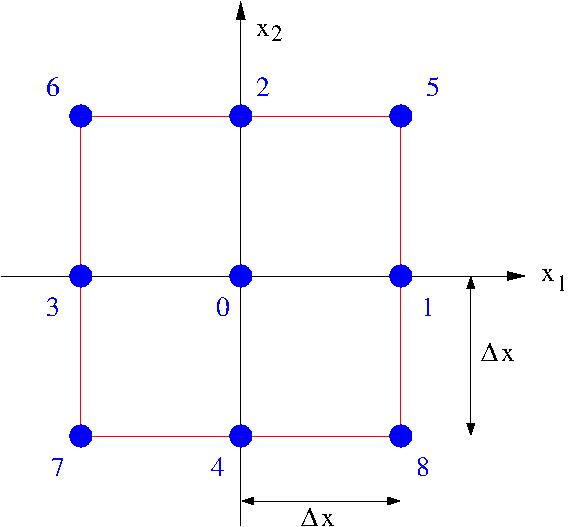}}
\caption{Particle distribution $ \, f_j \, $ for $ \, 0 \leq j \leq 8 \, $
of the D2Q9 lattice Boltzmann scheme.} 
\label{d2q9stencil}
\end{center} \end{figure}

The D2Q9 lattice Boltzmann schemes uses a set of discrete velocities described 
in  Fig.~\ref{d2q9stencil}. 
A density distribution $ \, f_j \,$ is associated to each velocity $ \, v_j  \equiv \lambda \, e_j $,  
where 
 $ \, \lambda = {{\Delta x}\over{\Delta t}} \, $ is the  fixed numerical lattice velocity.
From this particle distribution, we construct a vector  $ \, m \, $  of moments  $ \, m_k \,$ according to 
\moneq \label{orsay-moments} 
m = M \, f 
\monend
with an inversible fixed  matrix $ \, M \,$ usually \cite{LL00}  given by 
\moneq \label{matrice-M} 
 M  \, = \, \left (\begin {array}{ccccccccc}
\displaystyle   1 &   1 & 1 & 1 & 1 & 1 & 1 & 1 & 1\cr
\displaystyle   0 &   \lambda & 0 & -\lambda & 0 & \lambda & -\lambda & -\lambda  & \lambda\cr
\displaystyle   0 &   0 & \lambda & 0 & -\lambda & \lambda & \lambda & -\lambda & -\lambda\cr
\displaystyle -  4 \, \lambda^2 &   -\lambda^2  &  -\lambda^2 &  -\lambda^2 & -\lambda^2 &
 2 \,\lambda^2  &  2 \,\lambda^2 &  2 \,\lambda^2 &  2 \,\lambda^2 \cr  
\displaystyle   0 &    \lambda^2 & -  \lambda^2 &  \lambda^2 & -  \lambda^2 & 0 & 0 & 0 & 0\cr
\displaystyle   0 &   0 & 0 & 0 & 0 & 
\lambda^2 & -\lambda^2 & \lambda^2 & -\lambda^2    \cr  
\displaystyle   0 &   -2 \, \lambda^3 & 0 & 2  \, \lambda^3 & 0 & \lambda^3 &
 - \lambda^3 & - \lambda^3 &  \lambda^3 \cr
\displaystyle   0 &   0 & -2 \, \lambda^3 & 0 & 2 \, \lambda^3 &
  \lambda^3 &  \lambda^3 & - \lambda^3 & -  \lambda^3 \cr
\displaystyle   4 \, \lambda^4 &  -2 \, \lambda^4 &  -2  \, \lambda^4& -2  \, \lambda^4 &
 -2  \, \lambda^4 &  \lambda^4 & \lambda^4  &  \lambda^4 &   \lambda^4   \end{array}  \right) \, . 
\monend 
The three first moments for the  density and momentum are defined according to 
\moneq \label{orsay-moments-conserves} 
\rho \, = \,  \sum^8_{j=0} f_{j} \, = \, m_0 \,, \quad  
j_x \equiv \rho \, u_x \,  = \,   \sum^8_{j=0} \lambda \, e^1_{j} \, f_{j} \, = \, m_1 \,, \quad  
j_y \equiv \rho \, u_y \, = \,   \sum^8_{j=0} \lambda \, e^2_{j} \, f_{j} \, = \, m_2 \, ,  
\monend    
where  $\, e^\alpha_{j} \, $ are the cartesian components  of the vectors $e_{j}$ introduced previously. 
For fluid mechanics applications, a set of  ``conserved variables'' $ \, W \, $
is defined by 
\moneq \label{orsay-moments-conserves-2} 
W  \, =  \,  (\rho \,,\, j_x \,,\, j_y )  \, . 
\monend    
%

The moment distribution at equilibrium  $ \, m^{\rm eq} \, $ 
is a function of the conserved variables only.    
In particular, the equilibrium value of the fourth moment associated to kinetic energy 
(see the fourth line in the matrix $ \, M \, $ presented in (\ref{matrice-M})) 
is parametrized with the help of a constant~$  \alpha $. 
For the equilibrium of the last moment (of fourth order), we introduce 
 a constant $ \, \beta  $, as in our previous work \cite {DLT15}.
Note that the standard D2Q9 scheme \cite{LL00} uses $\alpha=-2$ and $\beta=1$.
The vector $ \, m^{\rm eq} \, $ of equilibrium moments is defined according to: 
\moneq \label{orsay-moments-equilibre} 
m^{\rm eq}  \,= \, \big( \rho \,,\, j_x \,,\, j_y \,, \, \alpha \, \lambda^2 \, \rho   \,,\, 0 \,,\, 0 \,,\, 
- \lambda^2 \, j_x \,, \, - \lambda^2  \, j_y  \,,\, \beta  \, \lambda^4 \, \rho 
\big)^{\displaystyle \rm t} \, . 
\monend    
Appling the inverse of  relation  (\ref{orsay-moments})  with the matrix $ \, M \, $ defined in (\ref{matrice-M}), we can 
explicit  herein all the components of  the vector $f^{\rm eq} $:
\begin{eqnarray}
\label{orsay-equilibre-fj-2} 
 {f^{\rm eq}  } \, : \,  \left\{ \begin{array} {rcl} 
f_0^{\rm eq} &=&\displaystyle  {{\rho} \over{9}}    \, \, \, \Big[ 1-\alpha +\beta \, \Big], \\  \vspace{-4 mm} \\ 
f_1^{\rm eq} &=&\displaystyle  {{\rho} \over{36}} \, \, \,\Big[ 4 -\alpha-2\, \beta   + {{12  \, u_x }\over{\lambda}}   \Big],\\  \vspace{-4 mm} \\
f_2^{\rm eq} &=&\displaystyle  {{\rho} \over{36}} \, \, \,\Big[ 4 -\alpha-2 \,\beta   + {{12  \, u_y }\over{\lambda}}   \Big], \\  \vspace{-4 mm}\\
f_3^{\rm eq} &=&\displaystyle  {{\rho} \over{36}} \, \, \,\Big[ 4 -\alpha-2 \, \beta   - {{12  \, u_x }\over{\lambda}}   \Big],\\  \vspace{-4 mm} \\
f_4^{\rm eq} &=&\displaystyle  {{\rho} \over{36}} \, \, \,\Big[ 4 -\alpha-2 \, \beta   - {{12  \, u_y }\over{\lambda}}   \Big],\\  \vspace{-4 mm} \\
f_5^{\rm eq} &=&\displaystyle  {{\rho} \over{36}} \, \, \,\Big[ 4 +2 \, \alpha+\beta   +{{3}\over{\lambda}}  \Big(u_x+u_y \Big)  \Big],
\\  \vspace{-4 mm} \\
f_6^{\rm eq} &=&\displaystyle  {{\rho} \over{36}} \, \, \,\Big[ 4 +2 \, \alpha+\beta   +{{3}\over{\lambda}}  \Big(-u_x+u_y \Big)   \Big],
\\  \vspace{-4 mm} \\
f_7^{\rm eq} &=&\displaystyle  {{\rho} \over{36}} \, \, \,\Big[ 4 +2 \, \alpha+\beta   +{{3}\over{\lambda}}  \Big(-u_x-u_y \Big) \Big],
\\  \vspace{-4 mm} \\
f_8^{\rm eq} &=&\displaystyle  {{\rho} \over{36}} \, \, \,\Big[ 4 +2 \, \alpha+\beta   +{{3}\over{\lambda}}  \Big(u_x-u_y \Big)   \Big].
 \end{array} \right.  
\end{eqnarray}

%

The lattice Boltzmann scheme is composed of two 
fundamental steps~: relaxation and advection.
During the relaxation step, the conserved variables  $ \, W \equiv ( \rho \,,\, J_x  \,,\, J_y ) \,$ are   not   modified; 
the nonconserved moments $\, m_3 \, $ to  $\, m_8 \, $ relax towards  an equilibrium value~: 
\moneq \nonumber  
m_k^{\rm eq} \, = \, \psi_k (W) \quad  {\rm for} \,\,  k \geq 3 \,,  
\monend 
%
where the $ \psi_k $ are linear functions of the conserved moments given in $(\ref{orsay-moments-equilibre}).$ 
This step depends upon relaxation rates 
$ \, s_k \, $ for $ \, k \geq 3  $:
\moneq \nonumber 
 m_k^* \, = \,  m_k \,+\, s_k \, \big(  m_k^{\rm eq} \,-\,  m_k \big)  \,, 
\monend   
where superscript $*$ denotes the moment $m_k$ after relaxation step. 
Now using the matrix $M^{-1}$ the relaxation step becomes in the $f$ space~:
\moneq \label{from-m-to-f}
 f^*_i(x, \, t) \, = \, \sum_\ell \, M^{-1} _{i \, \ell} \,  m_\ell^{*}.
\monend 
During the advection step $f_i(x_j)$ is ``transported''
from the node $x_j$ according to the discrete 
velocity $v_i$ to the node $x_j+v_i \Delta t.$
Thus the evolution of populations $f_i \, \ \ 0\leq i \leq 8,\, $ at internal node $x$ is described by~:
\moneq \label{LB-schema} 
f_i(x,t+\Delta t) \, = \, f^*_i(x - v_i\Delta t,t) \, .
\monend

\bigskip {\bf 2) \quad Bounce back boundary conditions for the D2Q9 scheme}  

\begin{figure}[htbp!]         
\begin{center}

\includegraphics[width=.55 \textwidth, angle=0]{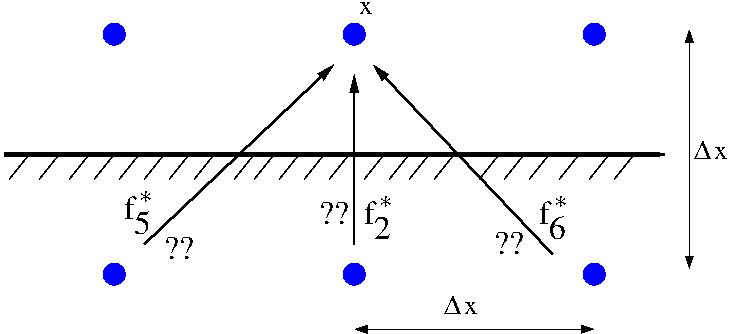}

\caption{Missing distribution functions to define the boundary scheme for the D2Q9 
on a boundary located at $y=0.$}
\label{d2q9bord}
\end{center} \end{figure}

Let us consider, without loss of generality,  the bottom boundary configuration 
as described in Fig.~\ref{d2q9bord}.
The values $f_i^{*}(x-v_i \Delta t)$  for $i \in \{2,5,6\}\equiv \mathcal{B}$  
to perform the scheme are unknown.

To impose a given velocity $(J_x,J_y)$ on the boundary we apply bounce back boundary condition~:   
\moneq \label{bb-classique} 
 \left \{ 
 \begin{array}{rcl}
  \displaystyle f_{2}(x,t+\Delta t) \,& = &
 \displaystyle f^*_{4}(x) \,+\, \frac{2}{3\lambda} J_y \Big(x,t+\frac{\Delta t}{2} \Big) \,,  \\  \vspace{-4 mm} \\
 \displaystyle f_{5}(x,t+\Delta t)  &= &
 \displaystyle  f^*_{7}(x) \,+\,  \frac{1}{6\lambda} ( J_x+ J_y)\Big(x-\frac{\Delta x}{2},t
+\frac{\Delta t}{2}\Big) \,,   \\  \vspace{-4 mm} \\
 \displaystyle f_{6}(x,t+\Delta t) &=& 
 \displaystyle f^*_{8}(x) \,+\,  \frac{1}{6\lambda} (-J_x+J_y) \Big(x+\frac{\Delta x}{2},t
+\frac{\Delta t}{2}\Big)  .  
\end{array}
\right.
\monend

\bigskip 
The bounce back scheme (\ref{bb-classique}) can be explained by a very  simple idea~: 
apply the internal scheme at the  boundary. If we focus on $f_2,$
we have 
 $ f_{2}(x, t+ \Delta t)  \,= \, f^*_{2}(x-(0,\, \Delta x), t ).$ 
For an internal node $x$ the particle distribution $ \, f_j \,$ is  close 
to the equilibrium idem for the particle distribution $ \, f_j^*\,$ after collision. 
So a  simple calculus leads to~:  
\moneq \nonumber 
f^*_{2}(x-(0,\, \Delta x))-f^*_{4}(x) \,=\, f_{2}^{\rm \rm eq}(x)-f_{4}^{\rm \rm eq}(x) +{\rm{O}}(\Delta x) \, . 
\monend  
Now we replace $f^{\rm{eq}}$ by their values given by $(\ref{orsay-equilibre-fj-2}),$ we get

\smallskip $ \displaystyle 
f^*_{2}(x-(0,\, \Delta x))-f^*_{4}(x) = 
 \, {{\rho} \over{36}} \, \,\Big[ 4 -\alpha-2 \,\beta   + {{12  \, u_y }\over{\lambda}}   \Big] \,
-\, {{\rho} \over{36}} \, \, \,\Big[ 4 -\alpha-2 \, \beta   
- {{12  \, u_y }\over{\lambda}}   \Big] + {\rm{O}}(\Delta x) $

\smallskip $ \displaystyle \qquad  \qquad  \qquad  \qquad \qquad  \quad    
= \, \frac{2}{3\lambda} J_y (x)  + {\rm{O}}(\Delta x) \, . $

\smallskip 
We remark here   
that $\rho u_y$ was substituted by the given function $J_y$ on the boundary which 
is the non-homogenous  
velocity to impose. An analogous calculus for the two other expressions 
gives the scheme (\ref{bb-classique}).

\monitem 
When the boundary condition is taken in a nonlinear way, the particle distribution at equilibrium 
is taken as a nonlinear function of the conserved quantities. Then taking appropriate
sums and differences of the corresponding relations that generalize  (\ref{orsay-equilibre-fj-2}), 
the extension of the bounce back conditions (\ref{bb-classique}) to the nonlinear framework is elementary. 
Our present analysis method is completely linear and we do not 
consider the nonlinear framework in this contribution. 

\begin{figure}[htbp!]         
\begin{center}
\includegraphics[width=.65 \textwidth, angle=0]{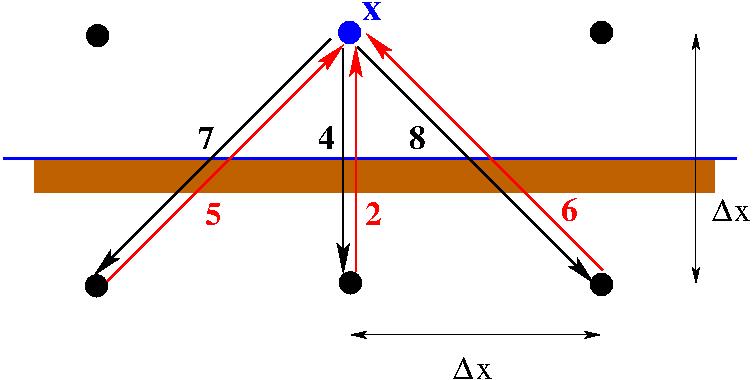}
\caption{Bounce back scheme for bottom boundary configuration with the D2Q9 scheme.}
\label{fig-bb}
\end{center}
\end{figure}

\monitem {\bf  Analysis of bounce back boundary condition.   }

In a previous work \cite{DLT15}  we made an analysis of the bounce back scheme 
using Taylor development 
which provides  a development of the velocity on the boundary node up to order two in space.

The main result is the following proposition. 

{\bf{Proposition 1.  

Expansion  of momentum at the node near the boundary up to order two}} 

\begin{figure}[htbp!]
\begin{center}
\centerline   {\includegraphics[height=.20  \textwidth] {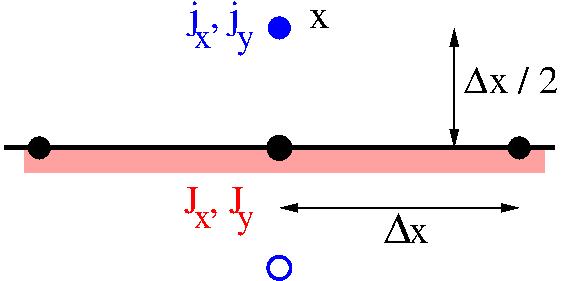}} 
\caption{Expansion of the momentum    $ \, (j_x ,\, j_y )  \,$ in the first cell
as a function of the data  $ \, ( J_x, \, J_y ) \, $ on the boundary.  }
\label{expansion} \end{center} \end{figure}

 The momentum   $ ( J_x, \, J_y ) $  (with upper-case letters)  is given on the boundary
(see Figure \ref{expansion}). In  this proposition, 
we  expand the momentum $ ( j_x, \, j_y ) $  (with lower-case letters) at the vertex $ \, x $,
located half a mesh size over the boundary,  in terms of these data. We have 
\moneq \label{dl-moments}
\left \{ \begin{array}{rcl}
  j_x&=&  \displaystyle J_x  
-  \frac{\Delta t}{2} \,  (4 \, \sigma_7+3) \,  \partial_t J_x 
+  \frac{\Delta x}{2} \,  \partial_y J_x 
\\  \vspace{-4 mm} \\    & & ~ \quad   
 \displaystyle +  \lambda \, \Delta x \,  \Big( 
 {{3 \alpha+2\beta+4}\over{6}}  \,  \sigma_7 - {{\alpha + 4}\over{6}} \, \Big( 2\, \sigma_7 + {3\over2} \Big) \Big) \, \partial_x \rho  
\\  \vspace{-4 mm} \\    & & ~ \quad        
+ \, \Delta x^2 \Big[  \alpha^0_{tt} \, \partial^2_t J_x + \alpha^0_{ty} \, \partial_t  \partial_y J_x 
+ \alpha^0_{xx} \, \partial^2_x J_x + \alpha^0_{yy} \, \partial^2_y J_x + \beta^0_{tx}  
\, \partial_t  \partial_x J_y   \,, 
\\  \vspace{-4 mm} \\    & & ~ \qquad \qquad 
 + \beta^0_{xy} \, \partial^2_{xy} J_y 
 + \gamma^0_{ty} \, \partial^2_{tx}  J_y 
+ \gamma^0_{xy} \, \partial_x \partial_y \rho \Big]
+  {\rm{O}} (\Delta x^3) \\  \vspace{-4 mm} \\ 
  j_y&=&\displaystyle  J_y 
-  \frac{\Delta t}{2} \, \partial_t J_y 
+  \frac{\Delta x}{2} \,   \partial_y J_y 
-  \frac{\Delta x}{12} \, ( \alpha + 4 ) \,  \partial_y \rho 
\\  \vspace{-4 mm} \\ 
       & & ~ \quad      + \,  \Delta x^2 \Big[ 
\theta^0_{tx} \, \partial_t \partial_x J_x + \theta^0_{xy} \, \partial_x \partial_y J_x 
+ \eta^0_{tt} \, \partial^2_t J_y +  \eta^0_{ty} \, \partial_t \partial_y J_y 
+ \eta^0_{xx} \, \partial^2_x J_y 
\\  \vspace{-4 mm} \\    & & ~ \qquad \qquad 
+ \eta^0_{yy} \, \partial^2_x J_y + \zeta^0_{ty} \, \partial_t \partial_y \rho + 
+ \zeta^0_{yy} \, \partial^2_y \rho \Big] +  {\rm{O}} (\Delta x^3)  \,. 
       \end{array}  \right. \monend
The coefficients that parametrize the second order terms in (\ref{dl-moments}) can be explicitly evaluated
and we have 
\moneq \label{bb0-coefs-2-bruts} \left\{   \begin{array}{l} \displaystyle 
\alpha^0_{tt} = {{1}\over{\lambda^2}} \, \Big( 6 \, \sigma_7^2+6 \, \sigma_7 + {13\over8} \Big) ,\, 
\alpha^0_{ty} = -{{1}\over{2 \, \lambda}} \, (2 \, \sigma_7+3 \, \sigma_4+3) ,  \\ \displaystyle 
\alpha^0_{xx} = {{1}\over24} \, (24 \, \sigma_4 \, \sigma_7
+8 \, \sigma_7 \, \sigma_8 + 12 \, \sigma_4 +8 \, \sigma_7 +15 ) ,\, 
\alpha^0_{yy} =  {{1}\over4} \, (2 \, \sigma_4+1) ,  \\ \displaystyle 
\qquad  \beta^0_{tx} = {{1}\over{12 \, \lambda}} \, (12 \, \sigma_4 \, \sigma_7 -4 \, \sigma_7 \, \sigma_8
+6 \, \sigma_4 -4 \, \sigma_7-9) , \,  \\ \displaystyle 
\beta^0_{xy} = -{{1}\over12} \, (12 \, \sigma_4 \, \sigma_7 -4 \, \sigma_7 \, \sigma_8
-4 \, \sigma_7-9) ,  \\ \displaystyle 
\qquad \gamma^0_{tx} = -{{1}\over12} \, \big(2 \, \alpha \, \sigma_3 \, \sigma_7
+6 \, \alpha \, \sigma_7^2 +12 \, \beta \, \sigma_7^2
+4 \, \beta \, \sigma_7 \, \sigma_8 -3 \, \alpha \, \sigma_3 -3 \, \alpha \, \sigma_7
\\ \displaystyle \qquad  \qquad \qquad \quad 
+6 \, \beta \, \sigma_7
-24 \, \sigma_7^2 -5 \, \alpha -\beta -40 \, \sigma_7-22 \big) , \\ \displaystyle 
\gamma^0_{xy} = {{\lambda}\over{36 \, (2 \, \sigma_7+1)}}  \, 
\big( 6 \, \alpha \, \sigma_3 \, \sigma_7^2 -6 \, \alpha \, \sigma_4 \, \sigma_7^2
+8 \, \alpha \, \sigma_7^2 \, \sigma_8 +4 \, \beta \, \sigma_3 \, \sigma_7^2 -12 \, \beta \, \sigma_4 \, \sigma_7^2
 \\ \displaystyle \qquad  \qquad  \qquad  \qquad 
+ 8 \, \beta \, \sigma_7^2 \, \sigma_8  
-9 \, \alpha \, \sigma_3 \, \sigma_7 -15 \, \alpha \, \sigma_4 \, \sigma_7
-4 \, \alpha \, \sigma_7^2 -2 \, \alpha \, \sigma_7 \, \sigma_8
 \\ \displaystyle \qquad  \qquad  \qquad  \qquad 
-6 \, \beta \, \sigma_3 \, \sigma_7 -6 \, \beta \, \sigma_4 \,  \sigma_7
+8 \, \sigma_3 \, \sigma_7^2 +24 \, \sigma_4 \, \sigma_7^2 -6 \, \alpha \, \sigma_4
 \\ \displaystyle \qquad  \qquad  \qquad  \qquad 
-29 \, \alpha \, \sigma_7 -6 \, \beta \, \sigma_7 -12 \, \sigma_3 \, \sigma_7 -36 \, \sigma_4 \, \sigma_7
-16 \, \sigma_7^2 -8 \, \sigma_7 \, \sigma_8 -9 \, \alpha
\\ \displaystyle \qquad  \qquad  \qquad  \qquad 
-24 \, \sigma_4 -92 \, \sigma_7 -36 \big) , \\ \displaystyle 
 \theta^0_{tx} = -{{1}\over{4 \, \lambda}} \, (2 \, \sigma_4 +1) , \, 
\theta^0_{xy} = {{1}\over4}  , \, 
\eta^0_{tt} = {{1}\over{8 \, \lambda^2}} , \,  
\eta^0_{ty} = -{{1}\over{6 \, \lambda}}  \, (\sigma_4+4) ,\, 
\eta^0_{xx} = {{4 \, \sigma_4 +1}\over24} , \\ \displaystyle 
\eta^0_{yy} = {{1}\over12}  \, (2 \, \sigma_4 +5)  ,\,\,\, 
\zeta^0_{ty} = {{1}\over24}  \,  (2 \, \alpha \, \sigma_3 +\alpha +8),\, \,\, 
\zeta^0_{xx} = -{{\lambda}\over24} \, (2 \, \sigma_4+1) \, (\alpha+4) ,  \\ \displaystyle 
\zeta^0_{yy} = -  {{\lambda}\over{72 \, (1 + 2 \, \sigma_7 ) }} 
\, \big( \, 6 \, \alpha \, \sigma_3 \, \sigma_7 -2 \, \alpha \, \sigma_4 \, \sigma_7
+4 \, \beta \, \sigma_3 \, \sigma_7 -4 \, \beta \, \sigma_4 \, \sigma_7 +2 \, \alpha \, \sigma_4
\\ \displaystyle \qquad  \qquad  \qquad  \qquad \quad 
+10 \, \alpha \, \sigma_7 
+8  \, \sigma_3 \, \sigma_7 +8 \, \sigma_4 \, \sigma_7 +5 \, \alpha +8 \, \sigma_4 +40 \, \sigma_7+20 \big) . 
 \end{array} \right. \monend 

%

\smallskip 
{\bf{Remark}}\\
 We  summarize here the demonstration established in our previous work~\cite{DLT15} 
 to use its results later and make this paper independent and clear.\\

\smallskip {\bf{Proof of Proposition 1.}}\\
 We write bounce back in general form~:
\moneq
 \displaystyle f^*_j(x,t+\Delta t)=f^*_\ell(x,t)+\xi_j(x',t'), \quad j\in \mathcal{B}  ,
 \monend
where $\ell$ is opposite of $j$ ({\it{i.e.}}  $v_j+v_\ell=0$) and $ \, \xi_j(x',t') \, $ 
is the given velocity on the boundary
as {\it e.g.} proposed  in (\ref{bb-classique}) to fix the ideas. 
If $j \notin \mathcal{B}$ the above equation is replaced by the internal scheme $(\ref{LB-schema}).$
%
By introducing  tables $T_{j,\ell}$,  and $U_{j,\ell}$ 
 the unified expression of the lattice Boltzmann scheme for a  node $x$ near the boundary is given by
\moneq  \label{bb-schema}
 f_j(x,t+\Delta t) =\sum_\ell T_{j,\ell} \, f_\ell^{*} (x,t) 
+ \sum_\ell U_{j,\ell} \, f^*_\ell(x-v_j\Delta t,t) + \xi_j 
\monend
where the matrix $U_{j,\ell}=1$ if $\ell=j \notin \mathcal{B}$ and $U_{j,\ell}=0$ if not. 
The matrix $\, U \,$ describes the ``internal'' numerical scheme (\ref{LB-schema})
whereas the matrix  $\, T \,$ takes into account the bounce back
boundary scheme (\ref{bb-classique}). 
For the particular bottom boundary with the D2Q9 scheme as presented in Fig.~\ref{fig-bb}, 
the tables $U$  and $T$ are~:
\smallskip \moneq \label{matrix-UT}
 \!\! \!\! \displaystyle U = \left( 
\begin{array}{ccccccccc}
1 & 0 & 0 & 0 & 0 & 0 & 0 & 0 & 0 \\ 
0 & 1 & 0 & 0 & 0 & 0 & 0 & 0 & 0 \\ 
0 & 0 & 0 & 0 & 0 & 0 & 0 & 0 & 0 \\ 
0 & 0 & 0 & 1  & 0 & 0 & 0 & 0 & 0 \\ 
0 & 0 & 0 & 0 & 1  & 0 & 0 & 0 & 0 \\ 
0 & 0 & 0 & 0 & 0 & 0 & 0 & 0 & 0 \\ 
0 & 0 & 0 & 0 & 0 & 0 & 0 & 0 & 0 \\ 
0 & 0 & 0 & 0 & 0 & 0 & 0 & 1  & 0 \\ 
0 & 0 & 0 & 0 & 0 & 0 & 0 & 0 &  1  
\end{array} \right) 
\quad \mbox{ and } \quad
 \displaystyle T = \left( 
\begin{array}{ccccccccc}
0 & 0 & 0 & 0 & 0 & 0 & 0 & 0 & 0 \\ 
0 & 0 & 0 & 0 & 0 & 0 & 0 & 0 & 0 \\ 
0 & 0 & 0 & 0 & 1 & 0 & 0 & 0 & 0 \\ 
0 & 0 & 0 & 0 & 0 & 0 & 0 & 0 & 0 \\ 
0 & 0 & 0 & 0 & 0 & 0 & 0 & 0 & 0 \\ 
0 & 0 & 0 & 0 & 0 & 0 & 0 & 1 & 0 \\ 
0 & 0 & 0 & 0 & 0 & 0 & 0 & 0 & 1 \\ 
0 & 0 & 0 & 0 & 0 & 0 & 0 & 0 & 0  \\ 
0 & 0 & 0 & 0 & 0 & 0 & 0 & 0 & 0 
\end{array}
\right).
\monend

\smallskip 
After linearization of the equilibrium, we can write the relaxation step as follows~: 
\moneq \nonumber 
m^* \,=\, J_0 \,\, m \,, 
\monend  
   with 
\moneq \label{matix-collision}
\displaystyle 
 J_0 = \left (
\begin{array}{ccccccccc}
\!\!\!\! 1 & \!\!\!\!\! 0 & \!\!\! 0 & \!\!\! 0 &  \!\!\! 0 &   \!\!\!\!\! 0 &  \!\!\!\!\! 0 &
\!\!\!\! 0\!\!\!\! &\!\!\!\! 0 \!\!\!\! \\ 
\!\!\!\!\!  0 & \!\!\!\!\! 1 &  \!\!\! 0 &  \!\!\! 0 &  \!\!\! 0 &   \!\!\!\!\! 0 &  \!\!\!\!\! 0 &
\!\!\!\! 0\!\!\!\! &\!\!\!\! 0 \!\!\!\! \\ 
\!\!\!\!\!  0 & \!\!\!\!\! 0 &  \!\!\! 1 &  \!\!\! 0 &  \!\!\! 0 &   \!\!\!\!\! 0 &  \!\!\!\!\! 0 &
\!\!\!\! 0\!\!\!\! &\!\!\!\! 0 \!\!\!\! \\ 
  \alpha \, s_3 \, \lambda^2 & \!\!\!\!\! 0 &  \!\!\! 0 &  \!\!\!\!\! 1-s_3 &  \!\!\!\!\! 0 &  
 \!\!\! 0 &  \!\!\! 0 & \!\!\!\! 0 \!\!\!\!&\!\!\!\! 0 \!\!\!\! \\ 
\!\!\!\!\!  0 & \!\!\! 0 &  \!\!\!\! 0 &  \!\!\! 0 &  \!\!\! 1-s_4  &   \!\!\!\!\! 0 &  \!\!\!\!\! 0 
&\!\!\!\! 0 \!\!\!\!&\!\!\!\! 0 \!\!\!\! \\ 
\!\!\!\!\!  0 & \!\!\! 0 &  \!\!\!\! 0 &  \!\!\! 0 &  \!\!\! 0  &    \!\!\!\!\! 1-s_4 &  \!\!\!\!\! 0 
&\!\!\!\! 0 \!\!\!\!&\!\!\!\! 0  \!\!\!\! \\ 
\!\!\!\!\!  0 & \!\!\!\!\! -s_7 \, \lambda^2  &  \!\!\! 0 &  \!\!\! 0 &  \!\!\!\!\! 0 &   \!\!\!\!\! 0 & 
\!\!\!\!1-s_7 &   \!\!\! 0 \!\!\!\!&\!\!\!\! 0 \!\!\!\! \\ 
\!\!\!\!\!  0 & \!\!\!\!\! 0 &  \!\!\! -s_7 \, \lambda^2  &   \!\!\! 0 &  \!\!\!\!\! 0 &   \!\!\!\!\! 0 &  
\!\!\! 0 & \!\!\!\!1-s_7 \!\!\!\!&\!\!\!\! 0 \!\!\!\! \\ 
\!\!\!\!  \beta \, s_8 \, \lambda^4 & \!\!\!\!\! 0 &  \!\!\! 0 &  \!\!\! 0 &  \!\!\! 0 &   
\!\!\!\!\! 0 &  \!\!\!\!\! 0 &\!\!\!\! 0\!\!\!\! &   1-s_8 
\end{array} \right).
\monend

\smallskip 
Then the unified LB scheme  (\ref{bb-schema}) for the bounce back becomes~: 
\moneq \label{eq-mom-bb}
 \displaystyle  
 \left\{ 
 \begin{array}{rcl} 
 m_k(x,t+\Delta t) &=&  (M T M^{-1} J_0)_{k,\ell} \ m_\ell (x,t) \\   \vspace{-4 mm} \\ 
&+&  (M_{k,\ell} U_{\ell,j} M^{-1}_{j,p} (J_0)_{p,q}) \ m_q (x-v_\ell \Delta t,t) + 
 M_{k,\ell} \, \, \xi_\ell  \, ,
\end{array}
\right. \monend
with an implicit summation on the repeated indices.
 We expand this relation to  order 0, 1 and 2. 

 \monitem {Analysis of bounce back at order zero}

Let us introduce the matrix \quad $  K \equiv I - M (T+U) M^{-1} J_0  $, 
  where $ \, I \, $ is the identity matrix. Then the equivalent equations 
for bounce back scheme  at order zero are the solution of 
\moneq \label{eqord0}
 \displaystyle   K m \,= \, M \, \xi   + {\rm{O}}(\Delta x) .
\monend 
For the D$2$Q$9$ scheme the matrix $K$ is given by~: 
\moneq \label{matrix-K}
K=\left( \begin{array}{ccccccccc}
0 & 0 & \frac{1}{\lambda}  & 0 & 0 & 0 & 0 & 0 & 0 \\ 
0 &   \frac{2-s_7}{3}  & 0 & 0 & 0 & 0 &  \frac{1 - s_7 }{3 \lambda^2 }  & 0 & 0  \\ 
0 & 0 & 1 & 0 & 0 & 0 & 0 & 0 & 0  \\ 
- s_3 \, \alpha \, \lambda^2   & 0 &   \tiny{\lambda ( 1 - s_7 )}  &    s_3 & 0 & 0 & 0 &    \frac{ 1 - s_7 }{\lambda} & 0  \\ 
0 & 0 &   -\frac{\lambda(1 + s_7)}{3}   & 0 &    s_4  & 0 &   0 &   {{1 - s_7}\over{3\, \lambda}} & 0 \\ 
0 &   \frac{\lambda( 2 - s_7) }{3}  & 0 & 0 & 0 &      s_4 &   \frac{1 - s_7}{3 \lambda}  & 0 & 0 \\ 
0 &   \frac{2\lambda^2  (1 + s_7 )}{3}    & 0 & 0 & 0 & 0 &  \frac{1+ 2 s_7 }{3} & 0 & 0  \\ 
0 & 0 & 0 & 0 & 0 & 0 & 0 & 1 & 0 \\ 
-\beta  s_8 \lambda^4 &  0 & - s_7 \lambda^3  & 0 & 0 & 0 & 0 &\lambda  ( 1 - s_7 )  &  s_8 
\end{array} \right) \, . 
\monend 
We remark that the matrix $K$ is  singular and the dimension of its kernel is equal to 1. 
In fact with $\mu_0 \, \equiv \,   \big( 1 ,\, 0 ,\, 0 ,\,  \alpha 
\lambda^2  ,\, 0  ,\, 0  ,\, 0  ,\,  0 ,\,  \beta \, \lambda^4  \big)^{\displaystyle \rm t}$,
we have  $ K \, \, \mu_0 = 0 $. 

In consequence, we have one compatibility  relation to satisfy,  
it is a linear combination of the equivalent equations of the internal scheme~: 
\moneqstar 
\lambda \, \big( \partial_t \rho + \partial_x j_x + \partial_y j_y \big) 
- \Big(\partial_t j_y + {{\alpha+4}\over{6}} \, \lambda^2 \,\partial_y  \rho  \Big) = 
{\rm O}(\Delta x) \, .  
\monendstar  
With given momenta $J_x$ and $J_y$  on the boundary, 
the density $ \, \rho \,$ remains still undefined by the boundary scheme. 
We develop the moments $m$ as~:
 $ m = m_0 + \Delta t \, m_1 + {\rm O}(\Delta t^2).$ We find that  the solution of (\ref{eqord0}) at order zero is~:
\moneqstar   
\displaystyle m_0=\left(\rho,J_x,J_y,\alpha \rho \lambda^2,0,0,-\lambda^2 J_x,-\lambda^2 J_y,\beta \rho \lambda^4
\right )^{\displaystyle \rm t}  \, . 
\monendstar  

\monitem {Analysis of bounce back at order one } 
 
Let introduce the matrix 
\moneqstar   
{  B_{k,p}^{\alpha} = 
\sum_{\ell,j,q} M_{k,\ell} \, U_{\ell,j}  \,v_j^{\alpha}  \,M_{j,q}^{-1}  \,(J_0)_{q,p} } , \quad   \alpha=1, \, 2 \, .
\monendstar  
Then the equivalent equations for bounce back scheme up to order one are solutions of 
\moneqstar 
  Km = M \, \xi  + \Delta t \left[ M \partial \xi - 
\partial_t m - B^{\alpha} \partial_\alpha m \right] 
+ {\rm{O}}(\Delta x^2) \, , 
\monendstar  
with $ \displaystyle m=m_0  + \Delta t \,  m_1 + {\rm{O}}(\Delta x^2)  $, and  
we have 
\moneqstar 
m_1= \Sigma \cdot \left(M \partial \xi 
- \partial_t m_0 - B^{\alpha} \partial_\alpha m_0 \right) \, . 
\monendstar  
The  matrix $\Sigma$ is given by~: 
\moneq
 \Sigma =  \left( \begin{array}{ccccccccc}
0 &  0 &   \!\!\!\!\!\!  \!\!\!\!\!\! 0  & 0 & 0 & 0 & 0 &  \!\!\!\!\!\!  \!\!\!\!\!\! 0 & 0 \\ 
0 &    2 + {{1}\over{s_7}}  &   \!\!\!\!\!\!  \!\!\!\!\!\! 0  & 0 &   0 & 0 & 
   {{1}\over{\lambda^2}}   \Big( 1 - {{1}\over{s_7}} \Big)  &  \!\!\!\!\!\!  \!\!\!\!\!\! 0 & 0 \\ 
0 &   0 &   \!\!\!\!\!\!  \!\!\!\!\!\!  1  & 0 & 0 & 0 & 0 & \!\!\!\!\!\!  \!\!\!\!\!\!  0 & 0 \\ 
0 &   0 &    \!\!\!\!\!\!  \!\!\!\!\!\!   {{\lambda}\over{s_3}}  \big( s_7 - 1 \big) 
  &   {{1}\over{s_3}}  & 0 & 0 & 0 & 
  \!\!\!\!\!\!  \!\!\!\!\!\!    {{1}\over{\lambda \, s_3}}  \big( s_7 - 1 \big) & 0 \\ 
0 &   0 &    \!\!\!\!\!\!  \!\!\!\!\!\!   {{\lambda}\over{3 \, s_4}}  \big( 1 + s_7  \big) 
  & 0 &    {{1}\over{s_4}}  & 0 & 0 & 
  \!\!\!\!\!\!  \!\!\!\!\!\!    {{s_7 - 1}\over{3 \, \lambda \, s_4}}  & 0 \\ 
0 &     - {{\lambda}\over{s_4}} 
 &   \!\!\!\!\!\!  \!\!\!\!\!\!  0  & 0 & 0 &    {{1}\over{s_4}}  & 0 &  \!\!\!\!\!\!  \!\!\!\!\!\! 0 & 0 \\ 
0 &      - 2  \lambda^2  \Big( 1 + {{1}\over{s_7}} \Big) 
 &   \!\!\!\!\!\!  \!\!\!\!\!\! 0  &  0 & 0 & 0 &     {{2}\over{s_7}}  - 1  & \!\!\!\!\!\!  \!\!\!\!\!\!  0 & 0 \\ 
0 &   0 &    \!\!\!\!\!\!  \!\!\!\!\!\! 0  & 0 & 0 & 0 & 0 & \!\!\!\!\!\!  \!\!\!\!\!\!  1 & 0 \\ 
0 &   0 &      \!\!\!\!\!\!  \!\!\!\!\!\!   \lambda^3  {{s_7}\over{s_8}}   & 0 & 0 & 0 & 0 & 
 \!\!\!\!\!\!  \!\!\!\!\!\!     {{\lambda}\over{s_8}}  \big( s_7 - 1 \big)  &     {{1}\over{s_8}}   
\end{array} \right).
\monend \\ 

  \monitem {Analysis of bounce back at order two}   

Let us introduce   
 the matrix 
$ \quad \displaystyle {  \widetilde{B}_{k,\ell}^{\alpha,\beta}
=\sum_{k,j,p,q} M_{k,p} \, U_{p,j} \,  v_j^{\alpha}  \, v_j^{\beta}  \, M_{j,q}^{-1}  \, (J_0)_{q,\ell} } $ , 
\quad   $1 \leq \alpha,\beta \leq 2.$ 

Then expand the various terms of the relation 
(\ref{eq-mom-bb}) up to second order.
We get the following equation~:
\moneq \label{ordre22} 
\left\{ \begin{array}{rcl}
 \displaystyle 
Km  = M \xi &+& \Delta t \,\, \left[ M \partial \xi - \partial_t m - B^{\alpha} \partial_\alpha m \right] \\
 \displaystyle 
&+& \frac{1}{2} \Delta t^2 \left[ M \partial^2\xi 
- \partial_t^2 m  \widetilde{B}^{\alpha,\beta} \partial_\alpha \partial_\beta m \right] + {\rm{O}} (\Delta t^3) \, . 
 \end {array}\right.  \monend

Now we develop $m$ as~: 
$ \displaystyle m=m_0  + \Delta t \,  m_1 + \Delta t^2 \,  m_2 + {\rm{O}}(\Delta x^3)  $. 
We get $ \,m_2 \,$ as solution of the equation $(\ref{ordre22})$~:
\moneq \nonumber  \displaystyle 
  m_2=K^{-1} \cdot \Big[ -\partial_t m_1 
- B^{\alpha} \partial_{\alpha} m_1  
+\frac{1}{2} \left(M\partial^2\xi -\partial^2_t m_0
+\widetilde{B}^{\alpha,\beta} \partial_\alpha \partial_\beta m_0\right ) \Big] .
\monend

\smallskip 
Finally the development of the velocities on the boundary nodes are given by the second 
and third components of $m=m_0  + \Delta t \,  m_1 + \Delta t^2 \,  m_2 + {\rm{O}}(\Delta x^3).$
When we explicit the conserved moments, the previous relation leads to the 
relation (\ref{dl-moments}). 
\hfill $\square$ 

\monitem 
We can use the  equivalent partial differential equations of the internal scheme
derived with the initial Taylor expansion method 
in order to express the relations (\ref{dl-moments}) and (\ref{bb0-coefs-2-bruts})
without the time derivatives. Recall that we have in the linearized case 
\moneq  \label{edp-2}   \left\{  \begin{array}{l}
\displaystyle \partial_t \rho  +    \partial_x j_x + \partial_y  j_y   = {\rm O}(\Delta x^2)  \,, \\  \vspace{-4 mm} \\ 
 \displaystyle   \partial_t j_x + c_0^2 \, \partial_x \rho 
- {{\lambda^2 \, \Delta t}\over3} \,\sigma_4 \, \big( \partial_x^2 +  \partial_y^2 \big) j_x
+ {{\lambda^2 \, \Delta t}\over6} \,\sigma_3 \, \alpha \, \partial_x \big(\partial_x j_x + \partial_y j_y ) 
 = {\rm O}(\Delta x^2)  \,, \\  \vspace{-4 mm} \\ 
 \displaystyle   \partial_t j_y + c_0^2 \, \partial_y \rho 
- {{\lambda^2 \, \Delta t}\over3} \,\sigma_4 \, \big( \partial_x^2 +  \partial_y^2 \big) j_y
+ {{\lambda^2 \, \Delta t}\over6} \,\sigma_3 \, \alpha \, \partial_y \big(\partial_x j_x + \partial_y j_y ) 
 = {\rm O}(\Delta x^2)  \,, 
\end{array} \right. \monend 
where $  \,  \displaystyle c_0=\lambda \sqrt{\frac{\alpha+4}{6}} \, $ is the speed of sound. 

In the result proposed in Proposition 1, we can take into account time dependent 
velocity  data on the boundary. In the following proposition, we transform
the first and second order time derivatives introduced in the expansion (\ref{dl-moments}) 
with the help of the partial differential equations presented at the relations (\ref{edp-2}).
We obtain after a tedious computation a new expression of the 
momentum in the cell directly close to the boundary. \\

{\bf{Proposition 2.  \\
A second expression of the velocity up to order two at the  boundary node}} 

We have the following expansions of the values $ \, j_x ,\, j_y \,$ at the boundary node 
in terms of the exact solution $\, J_x, \, J_y \,$ on the boundary (at $ \, x - {{\Delta x}\over{2}}$): 
\moneq \label{dl-moments-seconde}
\left \{ \begin{array}{rcl}
  j_x&=&  \displaystyle J_x  
+  \frac{\Delta x}{2} \,  \partial_y J_x 
+  {{\lambda \, \Delta x}\over6}  \, \big( 3 \alpha+2\beta+4 \big)  \,  \sigma_7  \, \partial_x \rho  
\\  \vspace{-4 mm} \\    & & ~  \quad        
+ \, \Delta x^2 \Big[  \widetilde{\alpha}^0_{xx} \, \partial^2_x J_x 
+ \widetilde{\alpha}^0_{yy} \, \partial^2_y J_x 
 + \widetilde{\beta}^0_{xy} \, \partial^2_{xy} J_y 
+ \widetilde{\gamma}^0_{xy} \, \partial_x \partial_y \rho \Big]
+  {\rm{O}} (\Delta x^3) \,,  \\  \vspace{-4 mm} \\ 
  j_y&=&\displaystyle  J_y 
+  \frac{\Delta x}{2} \,   \partial_y J_y 
-  \frac{\Delta x}{12} \, ( \alpha + 4 ) \,  \partial_y \rho 
\\  \vspace{-4 mm} \\ 
       & & ~ \quad      + \,  \Delta x^2 \Big[ \widetilde{\theta}^0_{xy} \partial_x \partial_y J_x 
+ \widetilde{\eta}^0_{xx} \partial^2_x J_y 
+ \widetilde{\eta}^0_{yy} \partial^2_x J_y 
+ \widetilde{\zeta}^0_{yy} \partial^2_y \rho \Big] +  {\rm{O}} (\Delta x^3) \,,
       \end{array}  \right. \monend
with 
\moneq \label{bb0-coefs-2-travailles} \left\{   \begin{array}{l} \displaystyle 
\widetilde{\alpha}^0_{xx} = {{1}\over48} \, (24 \, \alpha \, \sigma_3 \, \sigma_7
+72 \, \alpha \, \sigma_7^2 +48 \, \beta \, \sigma_7^2
+16 \, \beta \, \sigma_7 \, \sigma_8+36 \, \alpha \, \sigma_7+24 \, \beta \, \sigma_7 
\\ \displaystyle \qquad \qquad 
+16 \, \sigma_4 \, \sigma_7
 +96 \, \sigma_7^2+16 \, \sigma_7 \, \sigma_8-7 \, \alpha-4 \, \beta+48 \, \sigma_7 -6 ) ,\, 
\\ \displaystyle 
\widetilde{\alpha}^0_{yy} =  -{{1}\over12} \, (8 \, \sigma_4 \, \sigma_7-3) ,  \\ \displaystyle 
 \widetilde{\beta}^0_{xy} = {{1}\over48} \, (24 \, \alpha \, \sigma_3 \, \sigma_7+72 \, \alpha \, \sigma_7^2
+48 \, \beta \, \sigma_7^2 +16 \, \beta \, \sigma_7 \, \sigma_8+36 \, \alpha \, \sigma_7+24 \, \beta \, \sigma_7
-48 \, \sigma_4 \, \sigma_7
\\ \displaystyle \qquad \qquad 
+96 \, \sigma_7^2+16 \, \sigma_7 \, \sigma_8-7 \, \alpha-4 \, \beta+48 \, \sigma_7 ) , \, 
  \\ \displaystyle 
\widetilde{\gamma}^0_{xy} = {{\lambda}\over{72 \, (2 \, \sigma_7+1)}}  \, 
\big(12 \, \alpha \, \sigma_3 \, \sigma_7^2-36 \, \alpha \, \sigma_4 \, \sigma_7^2
 +24 \, \alpha \, \sigma_7^2 \, \sigma_8+8 \, \beta \, \sigma_3 \, \sigma_7^2-24 \, \beta \, \sigma_4 \, \sigma_7^2
\\ \displaystyle \qquad \qquad \qquad \qquad 
 +16 \, \beta \, \sigma_7^2 \, \sigma_8-18 \, \alpha \, \sigma_3 \, \sigma_7-18 \, \alpha \, \sigma_4 \, \sigma_7
 +24 \, \alpha \, \sigma_7^2-12 \, \beta \, \sigma_3 \, \sigma_7
\\ \displaystyle \qquad \qquad \qquad \qquad 
-12 \, \beta \, \sigma_4 \, \sigma_7+16 \, \sigma_3 \, \sigma_7^2
 -48 \, \sigma_4 \, \sigma_7^2+32 \, \sigma_7^2 \, \sigma_8+12 \, \alpha \, \sigma_7
-12 \, \beta \, \sigma_7
\\ \displaystyle \qquad \qquad \qquad \qquad 
-24 \, \sigma_3 \, \sigma_7-24 \, \sigma_4 \, \sigma_7+96 \, \sigma_7^2+9 \, \alpha+96 \, \sigma_7+36 \big) , \\ \displaystyle 
\widetilde{\theta}^0_{xy} = -{{\alpha}\over48}  , \, 
\widetilde{\eta}^0_{xx} = {{1}\over24}  ,\, 
\widetilde{\eta}^0_{yy} = -{{1}\over48}  \, (\alpha - 8)  ,\,  \\ \displaystyle 
\widetilde{\zeta}^0_{yy} = -  {{\lambda}\over{72 \, (1 + 2 \, \sigma_7 ) }} 
\, \big( 6 \, \alpha \, \sigma_3 \, \sigma_7-6 \, \alpha \, \sigma_4 \, \sigma_7
+4 \, \beta \, \sigma_3 \, \sigma_7-4 \, \beta \, \sigma_4 \, \sigma_7-6 \, \alpha \, \sigma_7
\\ \displaystyle \qquad \qquad \qquad \qquad \qquad 
+8 \, \sigma_3 \, \sigma_7 - 8 \, \sigma_4 \, \sigma_7-3 \, \alpha-24 \, \sigma_7-12 \big) . 
 \end{array} \right. \monend \\

\smallskip \smallskip  {\bf{Remark}} 

We have validated all the stationary coefficients of the equation $(\ref{dl-moments})$  
by different numerical test cases \cite{DLT15}.

\bigskip {\bf 3) \quad Towards a generalized first order bounce back boundary condition } 

To get a generalized  first order bounce back scheme the idea is 
to apply the internal scheme at the boundary 
\moneq
\left \{ \begin{array}{rcl}
 f_{5}(x, t+ \Delta t)  \,&=& \, f^*_{5}(x-(\Delta x,\, \Delta x), t ) ,\\ 
 f_{2}(x, t+ \Delta t)  \,&=&\, f^*_{2}(x-(0,\, \Delta x), t ) , \\
  f_{6}(x, t+ \Delta t)  \,&=& \, f^*_{6}(x+(\Delta x,\, -\Delta x), t ). 
  \end{array}
  \right.
 \monend
Then the expressions~:
\moneq \label{bb245768}
\left \{ \begin{array}{rlc}
  b_l  \,\, \,&\equiv&\, f^*_{5}(x-(\Delta x,\, \Delta x)) - f^*_{7}(x),\\ 
 b_m  \,&\equiv& \, f^*_{2}(x-(0,\, \Delta x))-f^*_{4}(x),  \\
 b_r \, \,\,&\equiv&\,  f^*_{6}(x+(\Delta x,\, -\Delta x), t ) - f^*_{8}(x)  , 
 \end{array}
 \right.
 \monend 
 are expanded at order one. 

\smallskip  {\bf{Remark}} \\ 
If  the expressions $b_l, b_m$ and $b_r$ are expanded only to order zero we get the traditional bounce back. 
 
\bigskip {\bf{Proposition 3.  First order bounce back}}  

For the boundary configuration described in Fig.~\ref{fig-bb}, the  
bounce back of first order scheme is given as follows~:  
\moneq \label{BB-O1} \left \{ \!\!  \begin{array}{rcl}
  f_{5}(x, t + \Delta t) \!\! \!&=& \!\!\!\! \displaystyle   \, f^*_{7}(x) \, +  
\,   {{1}\over{6 \, \lambda}} \,  \big( J_x +  J_y \big) \Big( x - {{\Delta x}\over{2}} \Big) 
+ {{1}\over{6}} \, \Big( q_x^* + q_y^* +  {{1}\over{\lambda}} \, \big(  j_x + j_y  \big)  \Big) 
\\  \vspace{-4 mm} \\ && \displaystyle  
\, + \, {{4 + \beta + 2 \, \alpha}\over{36}} \, 
\big( \rho(x) - \rho \big( x + (\Delta x, \Delta x) \big)  \big)  \,, 
 \\  \vspace{-4 mm} \\ 
f_{2}(x, t + \Delta t) \!  \!\! &=& \!\!\!\!  \displaystyle 
\, f^*_{4}(x)  \, +  \,   {{2}\over{3 \, \lambda}} \, J_y(x)   
\, -  {{1}\over{3}} \, \big( q_y^*(x) +  {{1}\over{\lambda}}\, j_y(x)  \big) 
\\  \vspace{-4 mm} \\ && \displaystyle  
\, + \, {{4 - 2 \, \beta - \alpha}\over{36}}  \, 
\big( \rho(x)  - \rho \big( x + (0,\, \Delta x) \big) \big)  \,, 
\\  \vspace{-4 mm} \\ 
f_{6}(x, t + \Delta t) \! \!\! &=& \!\!\!\!  \displaystyle \,  f^*_{8}(x) \, 
-  \, {{1}\over{6 \, \lambda}} \, \big( J_x -  J_y \big) \Big( x + {{\Delta x}\over{2}} \Big) 
\, +  {{1}\over{6}} \, \Big( \!\! -q_x^* + q_y^* \! + \! {{1}\over{\lambda}} \, 
\big( \! -j_x + j_y \big)  \Big) 
\\  \vspace{-4 mm} \\ && \displaystyle  
 \, + \,  \displaystyle {{4 + \beta + 2 \, \alpha}\over{36}}  \big( \rho(x) - \rho \big( x \! 
+   (-\Delta x, \Delta x) \big)  \big) .
\end{array} \right. \monend

\begin{figure}[htbp!]       
\begin{center}
\includegraphics[width=.55 \textwidth, angle=0]{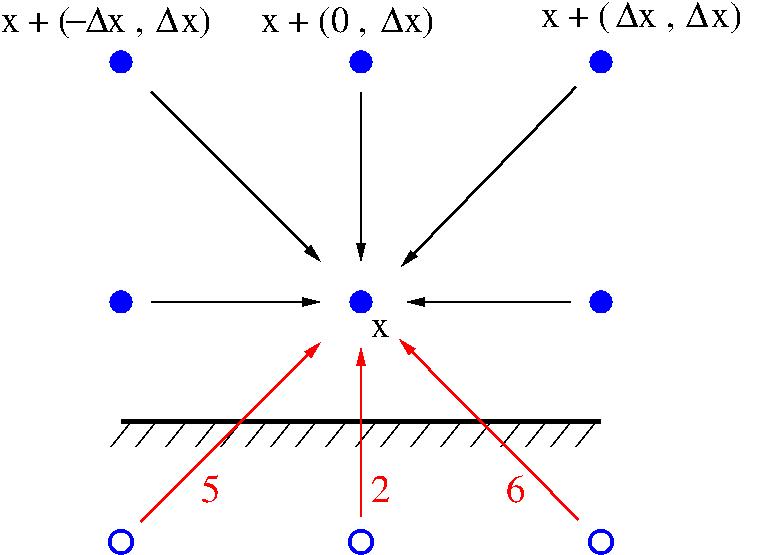}
\caption{Node near the boundary for the D2Q9 scheme.}
\label{grad_densite} \end{center} \end{figure}

\smallskip {\bf{Proof of Proposition 3.}}\\
We use Taylor expansion for the expressions $b_l, b_m$ and $b_r$ (see equations \ref{bb245768}), we get~: 
\moneq \label{BBdf} \left \{ \begin{array}{rcl}
b_l  \,\,\,  & = & \,  f^*_{5}(x) - f^*_{7}(x) 
+ {\rm d}f_5^{\rm \rm eq} (x) \smb (-\Delta x , \, -\Delta x ) + {\rm O}(\Delta x^2),  \\ \vspace{-4 mm} \\ 
 b_m  \,& = & \,  f^*_{2}(x)  - f^*_{4}(x)   + {\rm d}f_2^{\rm \rm eq} (x) \smb (0, \, -\Delta x ) 
+ {\rm O}(\Delta x^2), \\ \vspace{-4 mm} \\ 
b_r  \,\,\, &  = & \,  f^*_{6}(x) - f^*_{8}(x) 
+ {\rm d}f_6^{\rm \rm eq} (x) \smb (\Delta x , \, -\Delta x ) + {\rm O}(\Delta x^2) . 
\end{array}
\right.
\monend
With the help of (\ref{from-m-to-f}), we obtain the following exact expressions~:
\moneq \label{f256etoile} \left\{ \begin{array}{rcl}
  f^*_{5} - f^*_{7}  &\equiv&  \displaystyle {{1}\over{3 \, \lambda}}  \, \big( j_x + j_y \big) 
+ {{1}\over{6}} \, \big( q_x^*  + q_y^*  \big)\,\, ,   \\  \vspace{-4 mm} \\ 
 f^*_{2} - f^*_{4}  &\equiv & \displaystyle {{1}\over{3 \, \lambda}} \, j_y - {{1}\over{3}} \, q_y^* \,\,,   
 \\  \vspace{-4 mm} \\ f^*_{6} - f^*_{8} &\equiv & \displaystyle -{{1}\over{3 \, \lambda}}  
\, \big( j_x - j_y \big) - {{1}\over{6}} \, \big( q_x^*  - q_y^*  \big)\,\, . 
 \end{array} \right. \monend
On the other hand, we use the expressions $(\ref{orsay-equilibre-fj-2})$ for $f^{\rm{ eq}}_j,$ 
for $j$ equals $2$, $5$ and $6.$ 
So we have for D2Q9 particle equilibrium distribution~: 
\moneq \label{256eq} \left\{ \begin{array}{rcl}
 f_5^{\rm \rm eq}  \,&=&\, \displaystyle {{1}\over{36}} \, \big( 4 + 2 \, \alpha + \beta  \big) \, \rho +
 {{1}\over{12 \, \lambda}} \, \big( j_x + j_y \big),   \\  \vspace{-4 mm} \\
 f_2^{\rm \rm eq}  \,&=&\, \displaystyle {{1}\over{36}} \, \big( 4 - \alpha  - 2 \, \beta \big) \, \rho 
+ {{1}\over{3 \, \lambda}} \, j_y ,   \\  \vspace{-4 mm} \\
 f_6^{\rm \rm eq}  \,&=&\, \displaystyle {{1}\over{36}} \, \big( 4 + 2 \, \alpha + \beta \big) \, \rho
- {{1}\over{12 \, \lambda}} \, \big( j_x - j_y \big)  . 
 \end{array} \right. \monend 
In the expression $(\ref{BBdf})$ we need to expand ${\rm d}f_j^{\rm \rm eq}$  for $j$ equals $2$, $5$ and $6.$ 
Taking into account the above equations we need to develop the gradient 
of density and the gradient of momentum. \\
\monitem Gradient of density~:  using a Taylor expansion of the density around the node $x$
 (see Fig.~\ref{grad_densite}) we get~:
\moneq \label{dev_densite} \left \{ \begin{array}{rcl}
\, \nabla \rho \,\smb \, (\Delta x, \, \Delta x )  \, &\simeq& \, 
\rho \big( x + (\Delta x,\, \Delta x) \big) \,-\, \rho(x) \,+\,  {\rm O}(\Delta x^2)  \,, \\
\, \nabla \rho \, \smb  \, (0, \, \Delta x )  \, &\simeq& \, 
\rho \big( x + (0,\, \Delta x) \big) \,-\, \rho(x) \,+\,  {\rm O}(\Delta x^2) \,, \\
\, \nabla \rho \,  \smb \, (-\Delta x, \, \Delta x )  \, &\simeq& \, 
\rho \big( x + (-\Delta x,\, \Delta x) \big) \,-\, \rho(x) \,+\,  {\rm O}(\Delta x^2)   \,. 
\end{array} \right. \monend
%

\begin{figure}[htbp!]  
\begin{center}
\includegraphics[width=.45 \textwidth, angle=0]{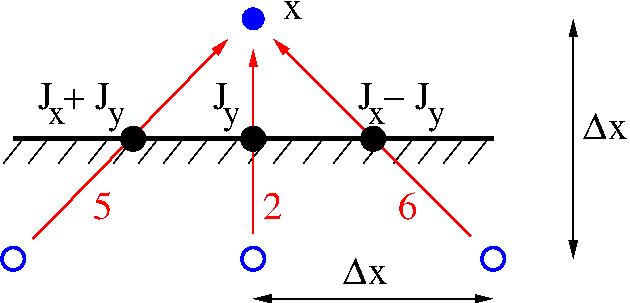}
\caption{Approximation of the momentum gradient
with finite differences.}
\label{grad_vitesse} \end{center} \end{figure}

\monitem Gradient of momentum~:  
using again Taylor expansion of the following combinations of the momentum 
around the node $x$ (see Fig.~\ref{grad_vitesse}) we get~:
\moneq \label{dev_vitesse} \!\!\!\! \left \{ \begin{array}{rcl}
 \nabla \big( j_x + j_y \big) \,  \smb \, (\Delta x , \, \Delta x ) \!\!  \! &\simeq& \!\!\! 
\displaystyle  2 \, \Big[ \big( j_x (x) + j_y(x) \big) 
- \big( J_x  +  J_y \big) \Big( x - {{\Delta x}\over{2}} \Big)  \Big] \,+\,  {\rm O}(\Delta x^2),  
 \\  \vspace{-4 mm} \\
\displaystyle  \nabla j_y \,  \smb \, (0, \, \Delta x )  \!\!\!   &\simeq& \!\!\! 
\displaystyle 2 \, \big( j_y(x) - J_y(x) \big) \,+\,  {\rm O}(\Delta x^2),   \\  \vspace{-4 mm} \\
 \displaystyle \nabla \big( j_x - j_y \big) \,  \smb \, (-\Delta x , \, \Delta x ) 
 \!\!  \!  &\simeq & \! \!\!  
\displaystyle 2 \, \Big[ \big( j_x (x) - j_y(x) \big)  - \big( J_x  -  J_y \big) 
\Big( x + {{\Delta x}\over{2}} \Big)  \Big] \,+\,  {\rm O}(\Delta x^2) .
\end{array} \right. \monend

\newpage 
\monitem Gradient of equilibrium particle distributions. 
Now using expressions (\ref{dev_densite}), (\ref{dev_vitesse})  and (\ref{256eq}) 
we obtain the following expressions for the gradient of equilibrium particle distributions~:
\moneq \label{dev_grad_f256} \!\!\!\! \left \{ \begin{array}{rcl} 
 \nabla f_5^{\rm eq}\smb (-\Delta x, \, -\Delta x )  \,&=&\, \displaystyle 
 {{1}\over{36}} \, \big( 4 + \beta + 2 \, \alpha \big) \, 
\big( \rho(x) - \rho \big( x + (\Delta x,\, \Delta x) \big)  \big) \\  \vspace{-4 mm} \\
&+& \displaystyle  {{1}\over{6 \, \lambda}} \,  \big[ \big( J_x  +  J_y \big) 
\big( x - {{\Delta x}\over{2}} \big) - 
\big( j_x (x) + j_y(x) \big)   \big] \,+\,  {\rm O}(\Delta x^2), \\  \vspace{-4 mm} \\
 \nabla f_2^{\rm eq}\smb (0, \, -\Delta x )  \,&=&\, \displaystyle 
 {{1}\over{36}} \, \big( 4 - 2 \, \beta - \alpha \big) \, 
\big( \rho(x)  - \rho \big( x + (0,\, \Delta x) \big) \big) \\  \vspace{-4 mm} \\
 &+ & \displaystyle {{2}\over{3 \, \lambda}} \,  \big( J_x (x) - j_y (x)\big)  \,+\,  {\rm O}(\Delta x^2)  , 
\\  \vspace{-4 mm} \\
\displaystyle  \nabla f_6^{\rm eq}\smb (\Delta x, \, -\Delta x )  \,&=&\, \displaystyle 
 {{1}\over{36}} \, \big( 4 + \beta + 2 \, \alpha \big) \, 
\big( \rho(x) - \rho \big( x + (-\Delta x,\, \Delta x) \big)  \big) \\  \vspace{-4 mm} \\
& -& \displaystyle  {{1}\over{6 \, \lambda}} \, 
 \big[ \big( J_x  -  J_y  \big) \big( x + {{\Delta x}\over{2}} \big) - 
\big( j_x (x) - j_y(x) \big)   \big] \,+\,  {\rm O}(\Delta x^2).
\end{array} \right. \monend
So by using the equations $(\ref{f256etoile})$ and $(\ref{dev_grad_f256})$ 
in the expressions given by $(\ref{BBdf})$ we get the following 
 expansion of the ``boundary gaps''~:  
%
%

\smallskip 
$ \displaystyle  b_l   \, \equiv   \,  f^*_{5}(x-(\Delta x,\, \Delta x)) -  f^*_{7}(x) $ 

\qquad  \qquad $ \displaystyle
\, =  {{1}\over{6 \, \lambda}} \, 
 \big( J_x +  J_y \big) \Big( x - {{\Delta x}\over{2}} \Big) 
 + \, {{1}\over{6}} \, \big( q_x^*(x) + q_y^*(x) 
+  {{1}\over{\lambda}} \, \big(  j_x(x) + j_y(x) \big)  \big) $ 

\qquad \qquad  \qquad $ \displaystyle
+ \, {{1}\over{36}} \big( 2 \, \alpha + \beta + 4 \big)    \, 
\big( \rho(x) - \rho \big( x + (\Delta x,\, \Delta x) \big)  \big)   + {\rm O}(\Delta x^2) \,, $

$ \displaystyle
b_m   \, \equiv     f^*_{2}(x-(0,\, \Delta x)) -  f^*_{4}(x)  $ 

\qquad  \qquad $ \displaystyle \,=   {{2}\over{3 \, \lambda}} \, J_y(x) \, - 
\, {{1}\over{3}} \, \big( q_y^*(x) +  {{1}\over{\lambda}}\, j_y(x)  \big)  $ 

\qquad \qquad  \qquad $ \displaystyle
+ \, {{1}\over{36}} \big( - \alpha - 2 \, \beta + 4 \big)    \, 
\big( \rho(x)  - \rho \big( x + (0,\, \Delta x) \big) \big)  + {\rm O}(\Delta x^2) \,, $

$ \displaystyle  
b_r  \, \equiv    
 f^*_{6}(x+(\Delta x,\, -\Delta x), t ) -  f^*_{8}(x)  $ 

\qquad  \qquad $ \displaystyle \, =  -{{1}\over{6 \, \lambda}} \, 
 \big( J_x  -  J_y \big) \Big( x + {{\Delta x}\over{2}} \Big) 
-   {{1}\over{6}} \, \big( q_x^*(x) - q_y^*(x) 
+  {{1}\over{\lambda}} \, \big(  j_x(x) - j_y(x) \big)  \big)  $ 

\qquad \qquad \qquad $ \displaystyle
+ \, {{1}\over{36}} \big( 2 \, \alpha + \beta + 4 \big)    \, 
\big( \rho(x) - \rho \big( x + (-\Delta x,\, \Delta x) \big)  \big)   + {\rm O}(\Delta x^2) $. 

\smallskip \noindent 
Thus the result of the proposition is a direct consequence of the above expressions. \hfill $\square.$

\bigskip {\bf 4) \quad Analysis of first order bounce back  } 

In this section we perform a formal analysis of  first order  bounce back scheme $(\ref{BB-O1})$  
given in the Proposition  $3.$
We consider here the same configuration as before ({\it{i.e.}} bottom boundary condition 
described in Fig.~\ref{d2q9bord}). 
In this case we give the 
 development of the velocity on the boundary node up to order two in space in the following proposition. 

\newpage 
\bigskip  \noindent 
{\bf{Proposition 4.  Expansion of the momentum for first order  bounce back}} 

We have the following expansion of momentum in the first cell up to order 2~: 
\moneq \label{dev_jxjy_o2} \left \{ \begin{array}{rcl}
 \displaystyle   j_x &=& J_x     \displaystyle +  \frac{\Delta x}{2} \,  \partial_y J_x
- 3 \, \Delta t \, \partial_t J_x
- \frac{\lambda \, \Delta x}{4} \, (4+\alpha) \,  \partial_y  \rho  \\  \vspace{-4 mm}  \\
 & & ~ + \, \Delta x^2 \, \Big[  \alpha^1_{tt} \, \partial^2_t J_x + \alpha^1_{ty} \, \partial_t  \partial_y J_x 
+ \alpha^1_{xx} \, \partial^2_x J_x + \alpha^1_{yy} \, \partial^2_y J_x + \beta^1_{tx}  \, \partial_t  \partial_x J_y 
\\  \vspace{-4 mm} \\ 
& &  \qquad  \qquad  + \beta^1_{xy}  \, \partial_x  \partial_y J_y 
+ \gamma^1_{tx} \, \partial_t \, \partial_x \rho + 
\gamma^1_{xy} \, \partial_x \, \partial_y \rho \Big] +  {\rm{O}} (\Delta x^3) \,, \\  \vspace{-4 mm}  \\
j_y &=& J_y  \displaystyle  + 
 {{\Delta x}\over2}  \, \partial_y J_y \, 
- \Delta t \, \partial_t J_y \, 
- {1\over6} \, \lambda \, \Delta x \,  (4+\alpha)  \, \partial_y \rho 
\\  \vspace{-4 mm} \\ \displaystyle
       & & ~ \quad      + \,  \Delta x^2 \Big[ 
\theta^1_{tx} \, \partial_t \partial_x J_x + \theta^1_{xy} \, \partial_x \partial_y J_x 
+ \eta^1_{tt} \, \partial^2_t J_y +  \eta^1_{ty} \, \partial_t \partial_y J_y 
+ \eta^1_{xx} \, \partial^2_x J_y 
\\  \vspace{-4 mm} \\    & & ~ \qquad \qquad 
+ \eta^1_{yy} \, \partial^2_x J_y + \zeta^1_{ty} \, \partial_t \partial_y \rho 
+ \zeta^1_{yy} \, \partial^2_y \rho \Big] +  {\rm{O}} (\Delta x^3) \, . 
\end{array}
\right.
\monend 
with 
\moneq \label{bb1-coefs-2-bruts} \left\{   \begin{array}{l} \displaystyle 
\alpha^1_{tt} =  {15\over{2 \, \lambda^2} } ,\,\,\, 
\alpha^1_{ty} = {{1}\over{4 \, \lambda}} \, (4 \, \sigma_7 - 6 \, \sigma_4 -11) ,\,\,\,   
\alpha^1_{xx} = {{1}\over24} \, (24 \, \sigma_4  + 4 \, \sigma_8 + 19 )  ,\, \\ \displaystyle 
\quad \alpha^1_{yy} =  {{1}\over4} \, (2 \, \sigma_4 + 1) ,\,  
\beta^1_{tx} = {{1}\over{12 \, \lambda}} \, (12 \, \sigma_4 \, -2 \, \sigma_8 - 11 ) , \,  \\ \displaystyle 
\beta^1_{xy} = {{1}\over12} \, (2 \, \sigma_8 - 6 \, \sigma_4 + 12 )   ,\, \,\,  
\gamma^1_{tx} = {{1}\over{12 \, \lambda}} \, (17 \, \alpha + \beta
-2 \, \beta \, \sigma_8 +2 \, \alpha \, \sigma_3+72)  ,\, \\ \displaystyle 
\quad \gamma^1_{xy} = {{\lambda}\over36} \, ( 88+28 \, \alpha+12 \, \sigma_4+6 \, \beta+4 \, \sigma_8
+3 \, \sigma_4 \, \alpha+12 \, \alpha \, \sigma_7+12 \, \beta \, \sigma_7 + \alpha \, \sigma_8)  ,\, 
\\ \displaystyle 
\theta^1_{tx} = - {{1}\over{12 \, \lambda}} \,(5-4 \, \sigma_7 +6 \, \sigma_4) ,\, \,\,  \,\,
\theta^1_{xy} =  {{1}\over4}  ,\, \,\, \,\, 
\eta^1_{ty} =  -{{1}\over{12 \, \lambda}} \, (11+2 \, \sigma_4)  ,\,  \\ \displaystyle 
\eta^1_{xx} =  {{1}\over24} \, (1 + 4 \, \sigma_4)  ,\, \,\,
\eta^1_{xy} =  {{1}\over12} \, (5 + 2 \, \sigma_4)  ,\, \,\,
\zeta^1_{ty} =  {{1}\over24} \, (2 \, \alpha \, \sigma_3+16+3 \, \alpha)  ,\,  \\ \displaystyle 
\zeta^1_{xx} =   -{{\lambda}\over36} \, (-4 \, \sigma_7 +16 + \beta+12 \, \sigma_4+5 \, \alpha
+ \alpha \, \sigma_7  +2 \, \beta \, \sigma_7+3 \, \sigma_4 \, \alpha) ,\, \\ \displaystyle 
\quad \zeta^1_{yy} =   -{{\lambda}\over72} \, (11+2 \, \sigma_4) \, (4+\alpha) \, . 
\end{array} \right. \monend

\smallskip {\bf{Proof of Proposition 4.}}\\
We use here exactly the same method as in the proof of the Proposition 1. We begin by 
writing the first order bounce back given by (\ref{BB-O1}) for  $ j\in \mathcal{B} = 
\{ 5, \, 2 , \, 6 \}  $ in the following form~:
\moneq
 \displaystyle 
f^*_j(x,t+\Delta t) = \sum_\ell T_{j,\ell}  \, f^*_\ell(x,t)+\xi_j(x',t),
\monend
where the transmission matrix  $ \, T $ is now~: 
\begin{center}
$ T=  \left(  \begin{array}{ccccccccc}
0 & 0 & 0 & 0 & 0 & 0 & 0 & 0 & 0 \\ 
0 & 0 & 0 & 0 & 0 & 0 & 0 & 0 & 0 \\ 
0 & 0 & 1/3 & 0 & {  2/3} & -2/3 & -2/3 & 2/3 & 2/3 \\ 
0 & 0 & 0 & 0 & 0 & 0 & 0 & 0 & 0 \\ 
0 & 0 & 0 & 0 & 0 & 0 & 0 & 0 & 0 \\ 
0 & -1/6 & -1/6 & 1/6 & 1/6 & 2/3 & 0 & {  1/3} & 0 \\ 
0 & 1/6 & -1/6 & -1/6 & 1/6 & 0 & 2/3 & 0 & {  1/3}  \\ 
0 & 0 & 0 & 0 & 0 & 0 & 0 & 0 & 0  \\ 
0 & 0 & 0 & 0 & 0 & 0 & 0 & 0 & 0 
\end{array}
\right)
$ 
\end{center}
and  $\xi_j(x', \, t)$ is the given data  for $ { j \in \mathcal{B} = \{ 5, \, 2 , \, 6 \} } $~: 
\moneq \label{dev_sigma}  \left \{ \begin{array}{rcl} 
 \displaystyle  \xi_5  \, &= & \, \displaystyle 
 {{1}\over{6}} \, \Big( 1 - {{\Delta x}\over{2}} \partial_x + {{\Delta x^2}\over{6}} \partial_x^2
\Big) \big( J_x +   J_y \big) \\\vspace{-4 mm}  \\
 && \qquad  + \, \displaystyle  {{4 + \beta + 2 \, \alpha}\over{36}} \,
\Big( - \Delta x ( \partial_x + \partial_y)  -  {{\Delta x^2}\over{2}} \big(  \partial_x + \partial_y) ^2 \Big) 
\rho  + {\rm O}(\Delta x^3), \\ \vspace{-4 mm}  \\
\displaystyle    \xi_2  \,& =&  \displaystyle  \,   {{2}\over{3 \, \lambda}} \,  J_y  
+  {{4 - 2 \, \beta - \alpha}\over{36}} \,
\Big( - \Delta x \, \partial_y -  {{\Delta x^2}\over{2}} \partial_x^2 \Big) 
\rho(x, \, t)  + {\rm O}(\Delta x^3) ,\\ \vspace{-4 mm}  \\
\xi_6  \,& =&  \, \displaystyle
 -{{1}\over{6}} \, \Big( 1 - {{\Delta x}\over{2}} \partial_x + {{\Delta x^2}\over{6}} \partial_x^2
\Big) \big( J_x -  J_y \big) \\ \vspace{-4 mm}  \\
 && \qquad + \, \displaystyle {{4 + \beta + 2 \, \alpha}\over{36}} \,
\Big( \Delta x ( \partial_x - \partial_y)  -  {{\Delta x^2}\over{2}} \big( ( \partial_x - \partial_y) ^2 \Big) 
\rho  + {\rm O}(\Delta x^3) .
\end{array} \right. \monend
Note here that the transmission matrix $T$ is modified to take into account the new given data $\xi$. 
Thus we can write the unified LB scheme (11) for the first order bounce back~:
\moneq \label{eq-mom-bb-bis}
 \displaystyle  
 \left\{ 
 \begin{array}{rcl} 
 m_k(x,t+\Delta t) &=&  (M T M^{-1} J_0)_{k,\ell} \ m_\ell (x,t) \\   \vspace{-4 mm} \\ 
&+&  (M_{k,\ell} U_{\ell,j} M^{-1}_{j,p} (J_0)_{p,q}) \ m_q (x-v_\ell \Delta t,t) + 
 M_{k,\ell} \xi_\ell  \, ,
\end{array} \right. \monend
where matrices $J_0$ and $U$ are given by $(\ref{matix-collision})$ 
and  $(\ref{matrix-UT})$ respectively.

Then as for the Proposition 1, we expand this relation for $m$ at order 0, 1 and 2 ~:
\moneq \nonumber  
m=m_0  + \Delta t \,  m_1 + \Delta t^2 \,  m_2 + {\rm{O}}(\Delta x^3) \,,
\monend  
and we expand also the equation $(\ref{eq-mom-bb-bis})$ at order 0, 1 and 2. 
At  order zero, we have to solve 
$ \displaystyle{   K  \, m_0 \,= \, M  \, \xi },$ where the matrix $K$ is  given by 
$ {   K \equiv I - M (T+U) M^{-1} J_0 } \, $~: 
\moneq \nonumber   
K=\left( \begin{array}{ccccccccc}
0 & 0 & \frac{1}{\lambda}  & 0 & 0 & 0 & 0 & 0 & 0 \\ 
0 &   \frac{1}{3}  & 0 & 0 & 0 & 0 &  0  & 0 & 0  \\ 
0 & 0 & 1 & 0 & 0 & 0 & 0 & 0 & 0  \\ 
- s_3 \, \alpha \, \lambda^2   & 0 &   0  &    s_3 & 0 & 0 & 0 &   0  & 0  \\ 
0 & 0 &   -\frac{2 \, \lambda}{3}   & 0 &    s_4  & 0 &   0 &   0  & 0 \\ 
0 &   \frac{\lambda}{3}  & 0 & 0 & 0 &      s_4 &  0   & 0 & 0 \\ 
0 &   \frac{\lambda^2  (1 + 3 \, s_7 )}{3}    & 0 & 0 & 0 & 0 &  s_7  & 0 & 0  \\ 
0 & 0 & \lambda^2  ( s_7 - 1)  & 0 & 0 & 0 & 0 & s_7 & 0 \\ 
-\beta  \, s_8  \, \lambda^4 &  0 & - \lambda^3  & 0 & 0 & 0 & 0 & 0  &  s_8 
\end{array} \right). 
\monend  
Note here that matrix $K$ here is different from that given by $(\ref{matrix-K})$ and is
still singular. 
So to solve the linear system $K m =g$, we must satisfy the compatibility condition which is
$ \displaystyle {   g_2 - \lambda \, g_0  \,= \, 0  } ,$ because the first and third lines 
of the matrix $K$ are proportional.
The linear space  $ \, \ker K  \, $ is generated by 
$ \displaystyle {    \mu  \, \equiv \,   \big( 1 ,\, 0 ,\, 0 ,\,  \alpha \lambda^2  ,\, 0  ,\, 0  ,\, 0  ,\,  0 ,\,  
\beta \, \lambda^4  \big)^{\displaystyle \rm t} } $. 
Thus the solutions of the equation  $K  \, m_0 \,= \, M \, \xi \, $  can be written as  
\quad $ \displaystyle {    
m_0 \,= \, \rho \, \mu  \,+\, \Sigma  \, M \, \xi }$ 
\qquad with 
\moneq \nonumber  
  { \Sigma = \left( \begin{array}{ccccccccc}
0 & 0 & 0 & 0 & 0 & 0 & 0 & 0 & 0 \\ 
0 &  3 & 0 & 0 & 0 & 0 &  0  & 0 & 0  \\ 
0 & 0 & 1 & 0 & 0 & 0 & 0 & 0 & 0  \\ 
0  & 0 &   0  &  1 /  s_3 & 0 & 0 & 0 &   0  & 0  \\ 
0 & 0 &   \frac{2 \, \lambda}{3 \, s_4}   & 0 &   1 / s_4   & 0  &   0 &   0  & 0 \\ 
0 &   - \lambda / s_4  & 0 & 0 & 0 &   1 / s_4 &  0   & 0 & 0 \\ 
0 &   - \frac{\lambda^2  (1 + 3 \, s_7 )}{s_7}    & 0 & 0 & 0 & 0 &  1 / s_7  & 0 & 0  \\ 
0 & 0 &  - \frac{\lambda^2  ( s_7 - 1) }{s_7}    & 0 & 0 & 0 & 0 & 1 / s_7  & 0 \\ 
0 & 0 & \frac{\lambda^3}{s_8} & 0 & 0 & 0 & 0 & 0  & 1 / s_8 
\end{array} \right).} 
\monend 
At  order zero,  we get~:  
\moneq \nonumber 
 {  m_0=\left(\rho, \, J_x, \, J_y, \,  \alpha \, \rho , \, 0, \, 0, \, -J_x / \lambda , \, - J_y/ \lambda,
 \,\beta  \, \rho \right )^{\displaystyle \rm t} }. 
\monend 
Note here that density $\rho$ is not fixed and $(J_x,J_y)$ are given functions on the boundary.\\
 Going further in the development of equation $(\ref{eq-mom-bb-bis}),$  we  get first order.
 Here the compatibility condition $ \,  g_2 - \lambda \, g_0  =  0  \,  $ is~: 
\moneq \nonumber  
\lambda \, \big( \partial_t \rho + \partial_x J_x + \partial_y J_y \big) 
- \Big(\partial_t J_y + {{\alpha+4}\over{6}} \, \lambda^2 \,\partial_y  \rho  \Big) = 
{\rm O}(\Delta x). 
\monend  
The above equation is a linear combination of the equivalent equations of the internal scheme. 
So this condition is satisfied. 
At  second  order, the compatibility condition still has the form~: 
\moneq \nonumber  
\lambda \, (\mbox{ conservation of the mass })  -  
 (\mbox{ conservation of momentum along} \,\, y)  = {\rm O}(\Delta x^2) \,. 
\monend   

Thus we can find the momentum development up to order two at the  boundary node 
at the vertex $x$ (inside the flow, located at $ \, \Delta x / 2  \,$ of the boundary). \hfill $\square$

\bigskip 
{\bf{Proposition 5.  Second expression of the momentum for first order bounce back }} 

When we replace in the expansion (\ref{dev_jxjy_o2}) the time derivatives by their
values obtained thanks to the partial equivalent equations (\ref{edp-2}), 
we obtain the following expansion of the value of momenta $ \, j_x \,$
and $ \,  j_y \,$ at the boundary node 
in terms of the exact solution $\, J_x, \, J_y \,$ on the boundary~: 
\moneq \label{dl-bb1-moments-seconde}
\left \{ \begin{array}{rcl}
  \!\!\!\!\!\!  j_x \!\!\!\!\! &=& \!\!\!\!\! J_x     \displaystyle 
\!+\!  \frac{\Delta x}{2} \,  \partial_y J_x  
 \!+\!  \Delta x^2  \Big[  \widetilde{\alpha}^1_{xx} \, \partial^2_x J_x 
\!+\!  \widetilde{\alpha}^1_{yy} \, \partial^2_y J_x \!+\!   \widetilde{\beta}^1_{xy}  \, \partial_x  \partial_y J_y 
\!+\!  \widetilde{\gamma}^1_{xy} \, \partial_x \, \partial_y \rho \Big] 
\!+\!  {\rm{O}} (\Delta x^3) 
\\  \vspace{-4 mm}  \\
  \!\!\!\!\!\!  j_y \!\!\!\!\! &=& \!\!\!\!\! J_y  \!+\!    \displaystyle 
 {{\Delta x}\over2}  \, \partial_y J_y   \!+\!  \Delta x^2 \Big[ 
 \widetilde{\theta}^1_{xy} \, \partial_x \partial_y J_x 
\!+\!  \widetilde{\eta}^1_{xx} \, \partial^2_x J_y \!+\!  \widetilde{\eta}^1_{yy} \, \partial^2_x J_y  
\!+\!  \widetilde{\zeta}^1_{yy} \, \partial^2_y \rho \Big] \!+\!   {\rm{O}} (\Delta x^3) \, . 
\end{array}
\right.
\monend

with 
\moneq \label{bb1-coefs-2-travailles} \left\{   \begin{array}{l} \displaystyle 
\widetilde{\alpha}^1_{xx} = {{1}\over24} \, 
(8 \, \alpha \, \sigma_3-5-2 \, \beta-4 \, \alpha+4 \, \sigma_8+4 \, \beta \, \sigma_8) ,\, 
\widetilde{\alpha}^1_{yy} =  -{{1}\over4} \, (2 \, \sigma_4 - 1 ) ,\,  \\ \displaystyle 
\qquad  \widetilde{\beta}^1_{xy} = {{1}\over12} \, (4 \, \alpha \, \sigma_3 
-2 \, \alpha+2 \, \sigma_8-6 \, \sigma_4  +2 \, \beta \, \sigma_8 - \beta - 1 ) ,\,   \\ \displaystyle 
\widetilde{\gamma}^1_{xy} = -{{\lambda}\over{6}}  \, (\beta +3 \, \alpha \, \sigma_7 
+2 \, \beta \, \sigma_7 +4 \, \sigma_7 + \alpha) ,\,
\widetilde{\theta}^1_{xy} = {{1}\over24}  \, (2 \, \alpha \, \sigma_3 -2 - \alpha) ,\,\\ \displaystyle 
\quad \widetilde{\eta}^1_{xx} = -{{1}\over24} \, (4 \, \sigma_4 -1) ,\, \quad 
\widetilde{\eta}^1_{yy} = {{1}\over24}  \, (2 \, \alpha \, \sigma_3-4 \, \sigma_4+2-\alpha)
,\,  \\ \displaystyle 
\widetilde{\zeta}^1_{xx} = -  {{\lambda}\over{72}} \, (8 \, \sigma_7+5 \, \alpha+2 \, \beta
+6 \, \alpha \, \sigma_7+4 \, \beta \, \sigma_7+12) . 
 \end{array} \right. \monend

The analysis of momentum at the boundary node obtained by first order bounce back scheme
proof shows that this scheme is more accurate than the simple bounce back scheme described by $(\ref{bb-classique}).$
In fact if we compare the analysis of the two schemes we see that 
for first order bounce back the order one terms are null (see equations $(\ref{dev_jxjy_o2})$.)
Moreover we note that with the following choice of the LB parameter  
$ \displaystyle  {\sigma_4 = {1\over4} }$ ({\it{i.e.}}  
$ {  - \frac{1}{4} \, ( 2 \, \sigma_4 - 1) \,=\, {1\over8} } $)   
the coefficients of  $\partial_y^2  J_x$ and $\partial_x^2  J_y$ 
in the equations $(\ref{dev_jxjy_o2})$ are null. Thus by this choice we get 
a quartic value  at the boundary for Poiseuille flow. 
This situation is not completely satisfactory and we propose in the following section 
to generalize the previous first order bounce back. 

 \newpage 
\bigskip {\bf 5) \quad Generalized  bounce back boundary scheme} 

Here we extend the first order bounce back scheme described by equations $(\ref{BB-O1})$ 
with the aim to cancel all the second 
order terms in the analysis of momentum at the boundary node given by $(\ref{dev_jxjy_o2}),$ allowing to get
a second order bounce back scheme. 
Let us introduce  
unknown parameters $a_k,$ $a_5,$  $a_6,$ $k_x$ and $k_y$ 
in the previous first order bounce back scheme (see equation $(\ref{BB-O1})$). 
Thus we get the following boundary scheme for bottom boundary (see Fig.~$\ref{fig-bb}$) ~: 
\moneq \label{BB-O2}
\left \{ \begin{array}{l} 
\displaystyle  f_{5}(x, t + \Delta t)  =   \displaystyle f^*_{7}  + {{1}\over{6 \, \lambda}} \, 
 \big( J_x +   J_y \big) \big( x - {{\Delta x}\over{2}}, \, t \big)
     \\  \vspace{-4 mm}  \\  \displaystyle \qquad \qquad 
+  {{{  a_5}}\over{6}} \, \big( q_x^* + q_y^*   + {{1}\over{\lambda}} \, 
\big(  j_x + j_y \big) \big) (x, \, t) 
+   {{{ k_5}}\over{36}} \,  \big( \rho(x, \, t) - \rho \big( x + (\Delta x, \Delta x) , \, t  \big)  \big)\,, 
 \\  \vspace{-4 mm}  \\ \displaystyle    f_{2}(x, t + \Delta t)   =   \displaystyle f^*_{4}(x)  
+   {{2}\over{3 \, \lambda}} \,  J_y  (x , \, t )
-  {{{  a_2}}\over{3}} \, \big( q_y^*   +  {{1}\over{\lambda}} \, j_y \big) (x, \, t) 
 \\    \displaystyle \qquad \qquad \qquad \qquad  
+   {{{ k_2}}\over{36}} \, {\big( \rho(x, \, t) - \rho \big( x + (0, \Delta x) , \, t  \big)  \big) } \,,
\\  \vspace{-4 mm}  \\  \displaystyle   f_{6}(x, t + \Delta t)   =   \displaystyle   f^*_{8}
 -  {{1}\over{6 \, \lambda}} \,  \big( J_x -   J_y \big) \big( x + {{\Delta x}\over{2}}, \, t \big)
 \\  \vspace{-4 mm} \\   \displaystyle  \quad 
 + {{{a_6}}\over{6}} \, \big( -q_x^* + q_y^*   + {{1}\over{\lambda}} \, 
\big(  -j_x + j_y \big) \big) (x, \, t) 
+   {{{ k_6}}\over{36}} \, \big( \rho(x, \, t) - \rho \big( x + (-\Delta x, \Delta x) , \, t  \big)  \big) .
\end{array}
\right.
\monend

\smallskip 
Note here that  we recover  simple bounce back scheme 
(essentially described by equations $(\ref{bb-classique})$) for all the parameters equal to zero 
({\it{i.e.}} $a_2=a_5=a_6=a_k=k_x=k_y=0$). 
If we choose all the parameters in the way proposed in relations (\ref{BB-O1}), {\it id est}
$ \, a_2=a_5=a_6=1 $, $ \, k_2 = 4-\alpha-2 \, \beta \, $ and 
$ \, k_5 = k_6 = 4+2 \, \alpha+\beta $, we 
recover the first order bounce back. 
It is possible to derive very long formal  expansions of the momenta
$ \, j_x \,$ and $ \, j_y \,$ in the first cell in terms of the boundary data, 
as in the relations (\ref{dl-moments}), (\ref{dl-moments-seconde}), (\ref{dev_jxjy_o2})
and (\ref{dl-bb1-moments-seconde}).

\monitem Analysis of the generalized bounce back 

For this extended bounce back we still have the following relation as for first order bounce back~: 
\moneq
\begin{array}{rcl}
\displaystyle  m_k(x,t+\Delta t) &=&  (M T M^{-1} J_0)_{k,\ell} \ m_\ell (x,t)  \\ 
\displaystyle  &+&  (M_{k,\ell} U_{\ell,j} M^{-1}_{j,p} (J_0)_{p,q}) \ m_q (x-v_\ell \Delta t,t) + 
 M_{k,\ell} \xi_\ell  \, . 
 \end{array}
 \monend 

 The matrix $U$ is unchanged but the transmission matrix $T$ satisfies~: 
\moneq \nonumber  
 { T=\left( \begin{array}{ccccccccc}
0 & 0 & 0 & 0 & 0 & 0 & 0 & 0 & 0 \\ 
0 & 0 & 0 & 0 & 0 & 0 & 0 & 0 & 0 \\ 
0 & 0 & {{a_2}\over{3}}  & 0 & {   1 - {{a_2}\over{3}}} & - {{2\, a_2}\over{3}} & - {{2\, a_2}\over{3}} 
 & {{2\, a_2}\over{3}}  & {{2\, a_2}\over{3}} \\ 
0 & 0 & 0 & 0 & 0 & 0 & 0 & 0 & 0 \\ 
0 & 0 & 0 & 0 & 0 & 0 & 0 & 0 & 0 \\ 
0 & -{{a_5}\over{6}} &  -{{a_5}\over{6}} & {{a_5}\over{6}}  & {{a_5}\over{6}} & 2 \, {{a_5}\over{3}} & 
0 & {  1 - {{2 \, a_5}\over{3}}} & 0 \\ 
0 & {{a_6}\over{6}} &  -{{a_6}\over{6}}  &  -{{a_6}\over{6}}  & {{a_6}\over{6}}  & 0 & {{2 \, a_6}\over{3}} &  
0 &  {  1 - {{2 \, a_6}\over{3}} }\\ 
0 & 0 & 0 & 0 & 0 & 0 & 0 & 0 & 0  \\ 
0 & 0 & 0 & 0 & 0 & 0 & 0 & 0 & 0 
\end{array} \right)}.
\monend  

 \newpage 
Now to get momentum development up to order two at the  boundary node we perform as in the previous section.
So  we expand this relation  at order 0, 1 and 2~:  
\moneq \nonumber  
m=m_0  + \Delta t \,  m_1 + \Delta t^2 \,  m_2 + {\rm{O}}(\Delta x^3) \,. 
\monend  
At order zero we introduce the matrix \quad $ {   K \equiv I - M (T+U) M^{-1} J_0 } \, .$
Remark that the matrix $\, K \,  $ depends on the  parameters  
$ \, a_2 $, $ \, a_5 $,  $ \, a_6 $,   $ \, a_k $,   $ \, k_x \,$ and  $ \, k_y $. 
Once again, we must solve the equation
\moneq \nonumber    
K  \, m_0 \,= \, M \, \xi   \,.  
\monend  
Its solution can be written as  
\moneq \nonumber  
m_0 \,= \, \rho \, \mu  \,+\, \Sigma \, \, M \, \xi \, .
\monend   
As in the previous  
case we get~: 
\moneq \nonumber   
{  m_0=\left(\rho, \, J_x, \, J_y, \,  \alpha \, \rho , \, 0, \, 0, \, 
-J_x / \lambda , \, - J_y/ \lambda,
 \,\beta  \, \rho \right )^{\displaystyle \rm t} }.
\monend   
 Note that in this case the compatibility condition at second order still has the form~:
 \begin{center} 
 $ \displaystyle \lambda \, $(conservation of the mass) $ - $ 
(conservation of momentum along $y$) $ = {\rm O}(\Delta x^2).$
\end{center} 
 The momentum $ \, \big( j_x \,,\, j_y \big) \,$ at the vertex $x$ located at $ \, \Delta x / 2  \,$ 
 can be expanded as powers of $ \, \Delta x$.  We fit the parameters  to get  no artefact at first order
 and to recover the Taylor expansion at the boundary. 

\bigskip \noindent 
{\bf{Proposition 6. \\ Interesting choice of the parameters of the generalized  bounce back}}  

If we take the parameters of the boundary scheme (\ref{BB-O2}) 
according to the relations
\moneq   \label{choix-parametres-ordre-1} \left \{ \begin{array}{rcl}
a_5   \,&=& \, a_6 , \\ 
 k_5 \, = \,  k_6 &=&\, \big (3 \, \alpha + 2\, \beta + 4 \big)  \, (1 - a_5 )  \, \sigma_7 
+  {1\over2} \,  \big (3 \, \alpha + 2\, \beta + 5 \big) \, a_5 +  {1\over2} \, (\alpha + 4)  \\
 k_2  \,&=&\,  2 \,  \big (3 \, \alpha + 2\, \beta + 4 \big)  \, (a_2 - 1) \,  \sigma_7 
- \big (3 \, \alpha + 2\, \beta + 4 \big) \, a_2 + 2 \, ( \alpha + 4) \,, 
\end{array} \right. \monend 
then we get the following expansion of the momentum $ \, \big( j_x \,,\, j_y \big) \,$ near the boundary~:
\moneq  \label{dev-bb2-basic} \left \{ \begin{array}{l} 
  j_x =  J_x     \displaystyle 
+  \frac{\Delta x}{2} \,  \partial_y J_x 
+   \displaystyle  ( 2 \, a_5 \,  \sigma_7 -2 \, \sigma_7 -2 -a_5 )  \, \Delta t \, \partial_t J_x 
\\  \vspace{-4 mm}  \\ \displaystyle \qquad \quad 
+ \frac{\lambda \, \Delta x}{6} \, (4+\alpha) \,  (2 \, a_5 \, \sigma_7-2 \, \sigma_7 -2-a_5) \, 
\partial_y  \rho  \\  \vspace{-4 mm}  \\ \displaystyle \qquad \quad 
   + \, \Delta x^2 \, \Big[  \alpha^2_{tt} \, \partial^2_t J_x  + \alpha^2_{ty} \, \partial_t \, \partial_x J_x
+ \alpha^2_{xx} \, \partial^2_x J_x + \alpha^2_{yy} \, \partial^2_y J_x + \beta^2_{tx}  \, \partial_t  \partial_x J_y 
\\  \vspace{-4 mm} \\ \displaystyle \qquad \qquad \qquad 
  + \beta^2_{xy}  \, \partial_x  \partial_y J_y 
+ \gamma^2_{tx} \, \partial_t \, \partial_x \rho + 
\gamma^2_{xy} \, \partial_x \, \partial_y \rho \Big] +  {\rm{O}} (\Delta x^3) \,, \\  \vspace{-4 mm}  \\
j_y =  J_y  \displaystyle  + 
 {{\Delta x}\over2}  \, \partial_y J_y \, 
- \Delta t \, \partial_t J_y \, 
- {1\over6} \, \lambda \, \Delta x \,  (4+\alpha)  \, \partial_y \rho 
\\  \vspace{-4 mm} \\ \displaystyle \qquad \quad 
   + \,  \Delta x^2 \Big[ 
\theta^2_{tx} \, \partial_t \partial_x J_x + \theta^2_{xy} \, \partial_x \partial_y J_x 
+ \eta^2_{tt} \, \partial^2_t J_y +  \eta^2_{ty} \, \partial_t \partial_y J_y 
+ \eta^2_{xx} \, \partial^2_x J_y 
\\  \vspace{-4 mm} \\  \displaystyle \qquad \qquad \qquad   
+ \eta^2_{yy} \, \partial^2_x J_y + \zeta^2_{ty} \, \partial_t \partial_y \rho 
+ \zeta^2_{yy} \, \partial^2_y \rho \Big] +  {\rm{O}} (\Delta x^3) \, . 
\end{array} \right. \monend 
with the following coefficients for the first component~: 
\moneq \nonumber  
\left\{   \begin{array}{l} \displaystyle 
\alpha^2_{tt} =   {{1}\over{\lambda^2}} \, ( 6 \, \sigma_7^2 +4 \, a_5^2 \, \sigma_7^2 -10 \, a_5 \, \sigma_7^2
+7 \, \sigma_7  -3 \, a_5 \, \sigma_7-4 \, a_5^2 \, \sigma_7+ {5\over2} +4 \, a_5 +a_5^2  ) ,\, 
\\ \displaystyle 
\,\,\,  \alpha^2_{ty} =  {{1}\over{4 \, \lambda}} \, 
(8 \, a_5 \, \sigma_7-4 \, \sigma_7-6 \, \sigma_4-7-4 \, a_5 )  ,\, 
\\ \displaystyle  
\quad \alpha^2_{xx} = -{{1}\over24} \, (8 \, a_5 \, \sigma_7 -8 \, \sigma_7 \, \sigma_8
+24 \, \sigma_4 \, a_5 \, \sigma_7 + 8 \, a_5 \, \sigma_7 \, \sigma_8 -8 \, \sigma_7 
-24 \, \sigma_4 \, \sigma_7 
\\ \displaystyle \qquad \qquad \qquad 
-12 \, \sigma_4 \, a_5 -4 \, a_5 - 4 \, a_5 \, \sigma_8- 15-12 \, \sigma_4 )  ,\, 
\qquad \alpha^2_{yy} =  {{1}\over4} \, (2 \, \sigma_4 + 1) ,\,  
\end{array}  \right.  \monend 
\moneq \nonumber 
\left\{   \begin{array}{l} \displaystyle 
\beta^2_{tx} = -{{1}\over{12 \, \lambda}} \, (9+12 \, \sigma_4 \, a_5 \, \sigma_7 +2 \, a_5+4 \, \sigma_7
-4 \, a_5 \, \sigma_7 +2 \, a_5 \, \sigma_8+4 \, \sigma_7 \, \sigma_8-6 \, \sigma_4
\\ \displaystyle \qquad \qquad \qquad 
-6 \, \sigma_4 \, a_5-12 \, \sigma_4 \, \sigma_7 -4 \, a_5 \, \sigma_7 \, \sigma_8 ) , \, 
\\ \displaystyle
\beta^2_{xy} = {{1}\over12} \, (12 \, \sigma_4 \, a_5 \, \sigma_7 +4 \, \sigma_7
-12 \, \sigma_4 \, \sigma_7 +4 \, \sigma_7 \, \sigma_8  -4 \, a_5 \, \sigma_7
-4 \, a_5 \, \sigma_7 \, \sigma_8+2 \, a_5
\\ \displaystyle \qquad \qquad 
-6 \, \sigma_4 \, a_5+2 \, a_5 \, \sigma_8+9)   ,\, 
\end{array}  \right.  \monend 
\moneq \nonumber 
\left\{   \begin{array}{l} \displaystyle 
\gamma^2_{tx} = {{1}\over12} \, \big( 30+7 \, \alpha + \beta +8 \, a_5^2+64 \, \sigma_7 
+3 \, \alpha \, \sigma_3 +34 \, a_5 -32 \, a_5 \, \sigma_7  -2 \, \beta \, a_5 \, \sigma_8
\\ \displaystyle \qquad \qquad 
-4 \, \beta \, \sigma_7 \, \sigma_8 -\alpha \, a_5 \, \sigma_3 -2 \, \alpha \, \sigma_7 \, \sigma_3  
-72 \, a_5 \, \sigma_7^2 +40 \, \sigma_7^2 +2 \, \alpha \, a_5 \, \sigma_7 \, \sigma_3  
\\ \displaystyle \qquad \qquad 
+4 \, \beta \, a_5 \, \sigma_7 \, \sigma_8 +32 \, a_5^2 \, \sigma_7^2 -32 \, a_5^2 \, \sigma_7
+2 \, \beta \, \sigma_7  +17 \, \alpha \, \sigma_7 -2 \, a_5 \, \sigma_7 \, \beta
\\ \displaystyle \qquad \qquad 
-9 \, \alpha \, a_5 \, \sigma_7  +8 \, a_5 \, \alpha+6 \, \sigma_7^2 \, \alpha  -4 \, \sigma_7^2 \, \beta
+2 \, a_5^2 \, \alpha +8 \, a_5^2 \, \sigma_7^2 \, \alpha -8 \, \alpha \, a_5^2 \, \sigma_7   
\\ \displaystyle \qquad \qquad 
-14 \, a_5 \, \sigma_7^2 \, \alpha+4 \, a_5 \, \sigma_7^2 \, \beta \big)  ,\, 
\\ \displaystyle 
\gamma^2_{xy} = -{{\lambda}\over36} \, \big( 60-6 \, \sigma_7 \, \alpha \, \sigma_4
+2 \, \sigma_7 \, \alpha \, \sigma_8+a_5 \, \alpha \, \sigma_8  +15 \, \alpha+8 \, \sigma_7 \, \sigma_8
+4 \, a_5 \, \sigma_8-12 \, \sigma_4 \, a_5
\\ \displaystyle \qquad \qquad 
-24 \, \sigma_4 \, \sigma_7  +24 \, \sigma_4+56 \, \sigma_7+28 \, a_5-56 \, a_5 \, \sigma_7
+6 \, \alpha \, \sigma_4-8 \, a_5 \, \sigma_7 \, \sigma_8 
\\ \displaystyle \qquad \qquad 
+24 \, \sigma_4 \, a_5 \, \sigma_7+24 \, \beta \, \sigma_7+38 \, \alpha \, \sigma_7-12 \, a_5 \, \sigma_7 \, \beta  -26 \, a_5 \, \sigma_7 \, \alpha+6 \, a_5 \, \beta+13 \, \alpha \, a_5 
\\ \displaystyle \qquad \qquad 
-2 \, \alpha \, a_5 \, \sigma_7  \, \sigma_8  +6 \, \alpha \, a_5 \, \sigma_7  \, \sigma_4
-3 \, a_5 \, \alpha \, \sigma_4 \big)  ,\, 
\end{array}  \right.  \monend 
and the coefficients for the second  component given according to~: 
\moneq \nonumber  
\left\{   \begin{array}{l} \displaystyle 
\theta^2_{tx} = - {{1}\over{4 \, \lambda}} \, 
{{(4 \, a_5 \, \sigma_7-3-2 \, a_5-6 \, \sigma_4) \, (2 \, \sigma_7 \, a_2-2 \, \sigma_7-a_2-1)}
\over{4 \, a_5 \, \sigma_7-2 \, a_5 +2 \, \sigma_7 \, a_2-6 \, \sigma_7-3-a_2}}  ,\, 
\\  \vspace{-4 mm}  \\  \displaystyle  
\theta^2_{xy} =  {{1}\over12} \, {1\over{4 \, a_5 \, \sigma_7-2 \, a_5 
+2 \, \sigma_7 \, a_2-6 \, \sigma_7-3-a_2}} \, \big( 12 \, \sigma_7 \, a_2 \, \sigma_4
-4 \, \sigma_7 \, a_2 \, \sigma_8+14 \, \sigma_7 \, a_2   
\\ \displaystyle \qquad \qquad \qquad \qquad
-6 \, \sigma_4 \, a_2+2 \, a_2 \, \sigma_8-7 \, a_2
+4 \, a_5 \, \sigma_7-18 \, \sigma_7+4 \, a_5 \, \sigma_7 \, \sigma_8  -12 \, a_5 \, \sigma_7 \, \sigma_4 
\\ \displaystyle \qquad \qquad \qquad \qquad
+6 \, a_5 \, \sigma_4-2 \, a_5 \, \sigma_8-2 \, a_5-9  \big) ,\, 
\end{array}  \right.  \monend 
\moneq \nonumber 
\left\{   \begin{array}{l} \displaystyle 
\eta^2_{ty} =  {{1}\over{12 \, \lambda}} \, {1\over{4 \, a_5 \, \sigma_7-2 \, a_5 
+2 \, \sigma_7 \, a_2-6 \, \sigma_7-3-a_2}} \, \big( -26 \, \sigma_7 \, a_2+4 \, \sigma_7 \, a_2 \, \sigma_8
+13 \, a_2
\\ \displaystyle \qquad \qquad \qquad \qquad
-2 \, a_2 \, \sigma_8  +12 \, \sigma_4 \, \sigma_7-4 \, a_5 \, \sigma_7 \, \sigma_8
-12 \, \sigma_4 \, a_5 \, \sigma_7-40 \, a_5 \, \sigma_7  +66 \, \sigma_7
\\ \displaystyle \qquad \qquad \qquad \qquad 
+20 \, a_5 +2 \, a_5 \, \sigma_8+6 \, \sigma_4 \, a_5+6 \, \sigma_4+33 \big)  ,\,  
\\ \displaystyle 
\eta^2_{xy} =  -{{1}\over12} \, {1\over{4 \, a_5 \, \sigma_7-2 \, a_5 
+2 \, \sigma_7 \, a_2-6 \, \sigma_7-3-a_2}} \, \big(-14 \, \sigma_7 \, a_2+4 \, \sigma_7 \, a_2 \, \sigma_8
+7 \, a_2
\\ \displaystyle \qquad \qquad \qquad \qquad 
-2 \, a_2 \, \sigma_8 -16 \, a_5 \, \sigma_7+12 \, \sigma_4 \, \sigma_7+30 \, \sigma_7
-12 \, \sigma_4 \, a_5 \, \sigma_7-4 \, a_5 \, \sigma_7 \, \sigma_8  
\\ \displaystyle \qquad \qquad \qquad \qquad 
+8 \, a_5+2 \, a_5 \, \sigma_8 +6 \, \sigma_4 \, a_5+6 \, \sigma_4+15 \big)  ,\,  
\end{array}  \right.  \monend 
\moneq \nonumber 
\left\{   \begin{array}{l} \displaystyle 
\zeta^2_{ty} =  {{1}\over24} \,  {1\over{4 \, a_5 \, \sigma_7-2 \, a_5 
+2 \, \sigma_7 \, a_2-6 \, \sigma_7-3-a_2}} \, \big(-48-9 \, \alpha-96 \, \sigma_7-6 \, \alpha \, \sigma_3
-24 \, a_2
\\ \displaystyle \qquad \qquad \qquad \qquad 
-24 \, a_5  +64 \, a_5 \, \sigma_7+32 \, \sigma_7 \, a_2+4 \, \beta \, a_5 \, \sigma_8
+2 \, \alpha \, a_5 \, \sigma_3-12 \, \alpha \, \sigma_7 \, \sigma_3  
\\ \displaystyle \qquad \qquad \qquad \qquad 
+32 \, \sigma_7^2 \, a_2
-32 \, a_5 \, \sigma_7^2-4 \, \alpha \, a_5 \, \sigma_7 \, \sigma_3  -8 \, \beta \, a_5 \, \sigma_7 \, \sigma_8
-18 \, \alpha \, \sigma_7+4 \, a_5 \, \sigma_7 \, \beta
\\ \displaystyle \qquad \qquad \qquad \qquad 
+18 \, a_5 \, \sigma_7 \, \alpha  
-6 \, a_2 \, \alpha-2 \, a_2 \, \beta+2 \, a_5 \, \beta-3 \, a_5 \, \alpha-4 \, \sigma_7 \, a_2 \, \beta
-8 \, \alpha \, a_2 \, \sigma_3  
\\ \displaystyle \qquad \qquad \qquad \qquad 
+24 \, \sigma_7^2 \, a_2 \, \alpha+16 \, \sigma_7^2 \, a_2 \, \beta
-4 \, a_2 \, \beta \, \sigma_8  -24 \, a_5 \, \sigma_7^2 \, \alpha-16 \, a_5 \, \sigma_7^2 \, \beta
\\ \displaystyle \qquad \qquad \qquad \qquad 
+8 \, \sigma_7 \, a_2 \, \beta \, \sigma_8  +16 \, \alpha \, \sigma_7 \, a_2 \, \sigma_3 \big)  ,\,  
\\ \displaystyle 
\zeta^2_{xx} =   {{\lambda}\over24} \,   {1\over{4 \, a_5 \, \sigma_7-2 \, a_5 
+2 \, \sigma_7 \, a_2-6 \, \sigma_7-3-a_2}} \, \big( 4 \, a_5 \, \sigma_7 \, \beta-16 \, \sigma_7
-8 \, \beta \, \sigma_7-12 \, \alpha \, \sigma_7  
\\ \displaystyle \qquad \qquad \qquad \qquad 
+10 \, a_5 \, \sigma_7 \, \alpha+24 \, a_5 \, \sigma_7-20-24 \, \sigma_4-5 \, a_5 \, \alpha-5 \, \alpha  
-6 \, \alpha \, \sigma_4-12 \, a_5-2 \, a_5 \, \beta \big) ,\, 
\\ \displaystyle 
\zeta^2_{yy} =   -{{\lambda}\over72} \,\,   {1\over{4 \, a_5 \, \sigma_7-2 \, a_5 
+2 \, \sigma_7 \, a_2-6 \, \sigma_7-3-a_2}} \, \big(-132-12 \, \sigma_7 \, \alpha \, \sigma_4
-2 \, a_5 \, \alpha \, \sigma_8  -33 \, \alpha
\\ \displaystyle \qquad \qquad \qquad \qquad 
+8 \, a_2 \, \sigma_8-8 \, a_5 \, \sigma_8
-24 \, \sigma_4 \, a_5-48 \, \sigma_4 \, \sigma_7-24 \, \sigma_4  -264 \, \sigma_7-40 \, a_2-92 \, a_5
\\ \displaystyle \qquad \qquad \qquad \qquad 
+208 \, a_5 \, \sigma_7+56 \, \sigma_7 \, a_2+48 \, \sigma_7^2 \, a_2  -48 \, a_5 \, \sigma_7^2
-6 \, \alpha \, \sigma_4+2 \, a_2 \, \alpha \, \sigma_8
\\ \displaystyle \qquad \qquad \qquad \qquad 
-16 \, \sigma_7 \, a_2 \, \sigma_8 
 +16 \, a_5 \, \sigma_7 \, \sigma_8+48 \, \sigma_4 \, a_5 \, \sigma_7-66 \, \alpha \, \sigma_7
+24 \, a_5 \, \sigma_7 \, \beta  +76 \, a_5 \, \sigma_7 \, \alpha
\\ \displaystyle \qquad \qquad \qquad \qquad 
-4 \, a_2 \, \alpha+6 \, a_2 \, \beta
-6 \, a_5 \, \beta-29 \, a_5 \, \alpha  -24 \, \sigma_7 \, a_2 \, \beta-10 \, \sigma_7 \, a_2 \, \alpha
+36 \, \sigma_7^2 \, a_2 \, \alpha
\\ \displaystyle \qquad \qquad \qquad \qquad 
+24 \, \sigma_7^2 \, a_2 \, \beta  -36 \, a_5 \, \sigma_7^2 \, \alpha
-24 \, a_5 \, \sigma_7^2 \, \beta-4 \, \sigma_7 \, a_2 \, \alpha \, \sigma_8 
 +4 \, a_5 \, \sigma_7 \, \alpha \, \sigma_8
\\ \displaystyle \qquad \qquad \qquad \qquad 
+12 \, a_5 \, \sigma_7 \, \alpha \, \sigma_4 -6 \, a_5 \, \alpha \, \sigma_4 \big) \, . 
\end{array} \right. \monend 

\bigskip  \noindent 
{\bf{Proposition 7. \\ A second expression of the momenta expansion for the generalized  bounce back}}  

If we take the parameters of the boundary scheme (\ref{BB-O2}) 
 according to (\ref{choix-parametres-ordre-1}), 
we can drop away the unstationary terms in (\ref{dev-bb2-basic}) with the help of the 
partial differential equations   (\ref{edp-2}).  
Then we get the following expansion of the momentum $ \, \big( j_x \,,\, j_y \big) \,$ near the boundary~:
\moneq \label{dev-bb2-final} \left \{ \begin{array}{rcl}
  \!\!\!\!\!\!  j_x \!\!\!\!\! &=& \!\!\!\!\! J_x     \displaystyle 
\!+\!  \frac{\Delta x}{2} \,  \partial_y J_x  
 \!+\!  \Delta x^2  \Big[  \widetilde{\alpha}^2_{xx} \, \partial^2_x J_x 
\!+\!  \widetilde{\alpha}^2_{yy} \, \partial^2_y J_x \!+\!   \widetilde{\beta}^2_{xy}  \, \partial_x  \partial_y J_y 
\!+\!  \widetilde{\gamma}^2_{xy} \, \partial_x \, \partial_y \rho \Big] 
\!+\!  {\rm{O}} (\Delta x^3) 
\\  \vspace{-4 mm}  \\
  \!\!\!\!\!\!  j_y \!\!\!\!\! &=& \!\!\!\!\! J_y  \!+\!    \displaystyle 
 {{\Delta x}\over2}  \, \partial_y J_y   
\\  \vspace{-4 mm}  \\ \displaystyle && \qquad 
+  \Delta x^2 \Big[ 
 \widetilde{\theta}^2_{xy} \, \partial_x \partial_y J_x 
+  \widetilde{\eta}^2_{xx} \, \partial^2_x J_y +  \widetilde{\eta}^2_{yy} \, \partial^2_x J_y  
+  \widetilde{\zeta}^2_{xx} \, \partial^2_x \rho 
+  \widetilde{\zeta}^2_{yy} \, \partial^2_y \rho \Big] +   {\rm{O}} (\Delta x^3) \, . 
\end{array}
\right.
\monend 
The associated  coefficients are given by the following relations~: 

\moneq \nonumber  
\left\{   \begin{array}{l} \displaystyle 
\widetilde{\alpha}^2_{xx} =   -{{1}\over24} \, (5-8 \, \sigma_7 \, \sigma_4-8 \, \sigma_7 \, \sigma_8
+6 \, \alpha \, \sigma_7 +4 \, \beta \, \sigma_7-6 \, a_5 \, \sigma_7 \, \alpha
+8 \, \beta \, a_5 \, \sigma_7 \, \sigma_8+12 \, \alpha \, a_5 \, \sigma_7 \, \sigma_3
\\ \displaystyle \qquad \qquad \quad 
+4 \, \alpha+2 \, \beta-4 \, a_5 \, \sigma_7 \, \beta-8 \, \beta \, \sigma_7 \, \sigma_8+4 \, \sigma_4
+8 \, \sigma_7+8 \, a_5 \, \sigma_7 \, \sigma_4+8 \, a_5 \, \sigma_7 \, \sigma_8
\\ \displaystyle \qquad \qquad \quad 
-2 \, \alpha \, \sigma_3
-8 \, a_5 \, \sigma_7+16 \, a_5 \, \sigma_7^2-16 \, \sigma_7^2-12 \, \sigma_7^2 \, \alpha-8 \, \sigma_7^2 \, \beta
-4 \, a_5 \, \sigma_8-4 \, a_5 \, \sigma_4
\\ \displaystyle \qquad \qquad \quad 
-4 \, \beta \, a_5 \, \sigma_8+8 \, a_5 \, \sigma_7^2 \, \beta
+12 \, a_5 \, \sigma_7^2 \, \alpha-12 \, \alpha \, \sigma_7 \, \sigma_3-6 \, \alpha \, a_5 \, \sigma_3 ) ,\, 
\\ \displaystyle 
\widetilde{\alpha}^2_{yy} =   {{1}\over12} \, ( 8 \, a_5 \, \sigma_7 \, \sigma_4 +3
-4 \, a_5 \, \sigma_4-2 \, \sigma_4-8 \, \sigma_7 \, \sigma_4 ) ,\, 
\end{array}  \right.  \monend 
\moneq \nonumber 
\left\{   \begin{array}{l} \displaystyle 
\widetilde{\beta}^2_{xy} =   -{{1}\over12} \, ( 1+12 \, \sigma_7 \, \sigma_4-4 \, \sigma_7 \, \sigma_8
+3 \, \alpha \, \sigma_7  +2 \, \beta \, \sigma_7-3 \, a_5 \, \sigma_7 \, \alpha
+4 \, \beta \, a_5 \, \sigma_7 \, \sigma_8+6 \, \alpha \, a_5 \, \sigma_7 \, \sigma_3
\\ \displaystyle \qquad \qquad \quad 
+2 \, \alpha +\beta-2 \, a_5 \, \sigma_7 \, \beta-4 \, \beta \, \sigma_7 \, \sigma_8+4 \, \sigma_7
-12 \, a_5 \, \sigma_7 \, \sigma_4+4 \, a_5 \, \sigma_7 \, \sigma_8 -\alpha \, \sigma_3
\\ \displaystyle \qquad \qquad \quad 
-4 \, a_5 \, \sigma_7
+8 \, a_5 \, \sigma_7^2-8 \, \sigma_7^2-6 \, \sigma_7^2 \, \alpha-4 \, \sigma_7^2 \, \beta-2 \, a_5 \, \sigma_8
+6 \, a_5 \, \sigma_4-2 \, \beta \, a_5 \, \sigma_8
\\ \displaystyle \qquad \qquad \quad 
+4 \, a_5 \, \sigma_7^2 \, \beta
+6 \, a_5 \, \sigma_7^2 \, \alpha -6 \, \alpha \, \sigma_7 \, \sigma_3-3 \, \alpha \, a_5 \, \sigma_3 ) ,\, 
\\ \displaystyle 
\widetilde{\gamma}^2_{xy} =   {{\lambda}\over6} \, ( 2 \, a_5 \, \sigma_7 \, \alpha+2 \, a_y5 \, \sigma_7 \, \beta
-a_5 \, \beta-5 \, \alpha \, \sigma_7 -4 \, \beta \, \sigma_7- a_5 \, \alpha-4 \, \sigma_7 )  ,\, 
\\ \displaystyle 
\widetilde{\theta}^2_{xy} =   -{{1}\over24} \, \,\,   {1\over{4 \, a_5 \, \sigma_7-2 \, a_5 
+2 \, \sigma_7 \, a_2-6 \, \sigma_7-3-a_2}} \, \big( -6-6 \, \alpha \, \sigma_7
+10 \, a_5 \, \sigma_7 \, \alpha-4 \, a_2 \, \alpha 
\\ \displaystyle \qquad \qquad \qquad \qquad 
-2 \, a_2 \, \beta-4 \, \sigma_7 \, a_2 \, \alpha-4 \, \sigma_7 \, a_2 \, \beta+2 \, a_5 \, \beta
+a_5 \, \alpha-4 \, \alpha \, a_2 \, \sigma_3  +24 \, \sigma_7^2 \, a_2 \, \alpha
\\ \displaystyle \qquad \qquad \qquad \qquad 
-4 \, a_2 \, \beta \, \sigma_8+16 \, \beta \, \sigma_7^2 \, a_2 
 +8 \, \alpha \, \sigma_7 \, a_2 \, \sigma_3+8 \, \sigma_7 \, a_2 \, \beta \, \sigma_8-4 \, a_2 \, \sigma_8
\\ \displaystyle \qquad \qquad \qquad \qquad 
-8 \, \beta \, a_5 \, \sigma_7 \, \sigma_8 
-20 \, \alpha \, a_5 \, \sigma_7 \, \sigma_3+12 \, \sigma_4 \, a_2
-3 \, \alpha+8 \, \sigma_7 \, a_2 \, \sigma_8   +4 \, a_5 \, \sigma_7 \, \beta
\\ \displaystyle \qquad \qquad \qquad \qquad 
-12 \, \sigma_7 -24 \, \sigma_7 \, a_2 \, \sigma_4+24 \, a_5 \, \sigma_7 \, \sigma_4  
-8 \, a_5 \, \sigma_7 \, \sigma_8
+6 \, \alpha \, \sigma_3-2 \, a_2-4 \, a_5
\\ \displaystyle \qquad \qquad \qquad \qquad 
+24 \, a_5 \, \sigma_7   -12 \, \sigma_7 \, a_2-32 \, a_5 \, \sigma_7^2
+32 \, \sigma_7^2 \, a_2+4 \, a_5 \, \sigma_8-12 \, a_5 \, \sigma_4 
\\ \displaystyle \qquad \qquad \qquad \qquad  
+4 \, \beta \, a_5 \, \sigma_8
-16 \, a_5 \, \sigma_7^2 \, \beta-24 \, a_5 \, \sigma_7^2 \, \alpha+12 \, \alpha \, \sigma_7 \, \sigma_3 
+10 \, \alpha \, a_5 \, \sigma_3 \big)   ,\, 
\end{array}  \right.  \monend 

\moneq \nonumber 
\left\{   \begin{array}{l} \displaystyle 
\widetilde{\eta}^2_{xx} =   {{1}\over24} \, \,   {1\over{4 \, a_5 \, \sigma_7-2 \, a_5 
+2 \, \sigma_7 \, a_2-6 \, \sigma_7-3-a_2}} \, \big( -6 \, \sigma_7-3 \, a_2-32 \, a_5 \, \sigma_7 \, \sigma_4
\\ \displaystyle \qquad \qquad \qquad \qquad  
+8 \, \sigma_7 \, a_2 \, \sigma_4  -3+12 \, \sigma_4+6 \, \sigma_7 \, a_2+16 \, a_5 \, \sigma_4
+24 \, \sigma_7 \, \sigma_4-4 \, \sigma_4 \, a_2 \big)  ,\, 
\\ \displaystyle 
\widetilde{\eta}^2_{yy} =   -{{1}\over24} \, \,   {1\over{4 \, a_5 \, \sigma_7-2 \, a_5 
+2 \, \sigma_7 \, a_2-6 \, \sigma_7-3-a_2}} \, \big( 6-24 \, \sigma_7 \, \sigma_4-6 \, \alpha \, \sigma_7
+10 \, a_5 \, \sigma_7 \, \alpha
\\ \displaystyle \qquad \qquad \qquad \qquad  
  -4 \, a_2 \, \alpha-2 \, a_2 \, \beta-4 \, \sigma_7 \, a_2 \, \alpha
-4 \, \sigma_7 \, a_2 \, \beta+2 \, a_5 \, \beta  +a_5 \, \alpha-4 \, \alpha \, a_2 \, \sigma_3
\\ \displaystyle \qquad \qquad \qquad \qquad  
+24 \, \sigma_7^2 \, a_2 \, \alpha-4 \, a_2 \, \beta \, \sigma_8+16 \, \beta \, \sigma_7^2 \, a_2 
  +8 \, \alpha \, \sigma_7 \, a_2 \, \sigma_3+8 \, \sigma_7 \, a_2 \, \beta \, \sigma_8
\\ \displaystyle \qquad \qquad \qquad \qquad  
-4 \, a_2 \, \sigma_8
  -8 \, \beta \, a_5 \, \sigma_7 \, \sigma_8-20 \, \alpha \, a_5 \, \sigma_7 \, \sigma_3-8 \, \sigma_4 \, a_2
 -3 \, \alpha+8 \, \sigma_7 \, a_2 \, \sigma_8
\\ \displaystyle \qquad \qquad \qquad \qquad  
+4 \, a_5 \, \sigma_7 \, \beta-12 \, \sigma_4+12 \, \sigma_7
+16 \, \sigma_7 \, a_2 \, \sigma_4   +8 \, a_5 \, \sigma_7 \, \sigma_4-8 \, a_5 \, \sigma_7 \, \sigma_8
\\ \displaystyle \qquad \qquad \qquad \qquad  
+6 \, \alpha \, \sigma_3   -2 \, a_2+8 \, a_5-12 \, \sigma_7 \, a_2-32 \, a_5 \, \sigma_7^2+32 \, \sigma_7^2 \, a_2
+4 \, a_5 \, \sigma_8 
\\ \displaystyle \qquad \qquad \qquad \qquad  
 -4 \, a_5 \, \sigma_4+4 \, \beta \, a_5 \, \sigma_8-16 \, a_5 \, \sigma_7^2 \, \beta
-24 \, a_5 \, \sigma_7^2 \, \alpha   +12 \, \alpha \, \sigma_7 \, \sigma_3+10 \, \alpha \, a_5 \, \sigma_3 \big) ,\, 
\\ \displaystyle 
\widetilde{\zeta}^2_{xx} =   {{\lambda}\over24} \, \,   
{{2 \, \sigma_7 \, a_2-2 \, \sigma_7-a_2-1 }\over{4 \, a_5 \, \sigma_7-2 \, a_5 
+2 \, \sigma_7 \, a_2-6 \, \sigma_7-3-a_2}} \, \big(-8 \, \beta \, \sigma_7-12 \, \alpha \, \sigma_7
-16 \, \sigma_7+8 \, a_5 \, \sigma_7
\\ \displaystyle \qquad \qquad \qquad \qquad  
+6 \, a_5 \, \sigma_7 \, \alpha+4 \, a_5 \, \sigma_7 \, \beta
-3 \, a_5 \, \alpha-2 \, \alpha-4 \, a_5-2 \, a_5 \, \beta-8 \big)   ,\, 
\\ \displaystyle  
\widetilde{\zeta}^2_{yy} =   -{{\lambda}\over24} \, \,   
{{ ( 2 \, \sigma_7-1 )^2 }\over{4 \, a_5 \, \sigma_7-2 \, a_5 
+2 \, \sigma_7 \, a_2-6 \, \sigma_7-3-a_2}} \, \big(3 \, \alpha+2 \, \beta+4 \big) \, \big( a_2-a_5 \big) 
\, . 
\end{array} \right. \monend

\bigskip {\bf 6) \quad Numerical test}  

\monitem
Let first consider Poiseuille flow driven by  a pressure gradient in the domain $\Omega=[1, N_x] \times [1,N_y].$  
So we apply ``anti-bounce back'' boundary condition at inlet $(i=1)$  and outlet $(i=N_x)$ 
of the channel to
impose pressure $\delta p$ and $-\delta p$ and simple bounce back on the bottom $j=1$ 
and the top $j=N_{y}$ of the domain to impose $u_{y}=0$. 
So the solid wall ($j_x=0$) of the Poiseuille solution is exactly at $\frac{\Delta x}{2}$ 
for the following condition~:
\moneqstar 
\sigma_{4} \, \sigma_{7} \, = \, -\frac{3}{8}\frac{\alpha+4}{\alpha+2 \beta-4} \,,
\monendstar  
as proposed in our previous contribution \cite{DLT10}.\\
Now we use the extended bounce back scheme described by $(\ref{BB-O2})$ instead of
the classical bounce back to impose the homogeneous Dirichlet boundary condition $u_{y}=0$.
In this case we find that the parabolic solution of the Poiseuille flow is null exactly at $\frac{\Delta x}{2}$
({\it{i.e.}} the solid wall $j_x=0$ is located exactly at $\frac{\Delta x}{2}$
below the first mesh vertex) for the following choice of the parameters  $ \, a_5 $, $ \, s_4 \, $ and $ \, s_7$~: 
\moneq \label{bb2-quartic} 
\sigma_{4} \, \sigma_{7} \, = \, {1\over16} \, {{ 8 \, a_5 \, \sigma_4 + 4 \, \sigma_4 -3}\over{ (a_5 - 1) }}
  \, . 
\monend  
This is because when the coefficient $ \, \widetilde{\alpha}^2_{yy} \,$ 
introduced in Proposition~7 is equal to  $\, {1\over8} $, we have the relation (\ref{bb2-quartic}). 
%
This proves numerically that this extended bounce back scheme is exact 
for Poiseuille flow test case. 
All the extra order terms of the expressions $ \, (\ref{dev-bb2-final}) \, $ are null and the developments
 of the $j_x$ momentum on the boundary node becomes~:
\moneq \nonumber
 \displaystyle        
j_x = J_x  +  \frac{\Delta x}{2} \,  \partial_y J_x+ \frac{\Delta x^2}{8} \,  \partial_y^2  J_x.  \\
\monend

\monitem
We consider  now the ``accordion'' test case introduced in our contribution \cite{DLT15}. 
In the rectangular domain  $\Omega = ]0,L[ \times ] 0,h [$, 
we introduce  periodic boundary conditions at $x=0$ and $x=L$.
For the boundaries at $y=0$ and $y=h$,  we impose 
$\, J_x(x,0)=J_x(x,h) = J_0\ \cos(2 k \pi \frac{x}{L}) \, $ and   
$\, J_y(x,0) = J_y (x,h) =0, $ for  $ \,  0 < x < L \,$ 
with  the integer $ \, k \, $ equal to 1 in our simulations. 
In the low velocity  regime  the
steady state is solution of the  Stokes equations 
\moneq \label{stokes} \displaystyle         
{\rm div} \, J = 0 \,, \,\,  - \nu \, \Delta J + \nabla  p = 0   \,. 
\monend
An analytic solution is  given by the following expressions. 
Introduce the function $ \, f(y) \, $ defined by 
\moneq \nonumber \displaystyle 
f(y) = \left\{  \begin {array} {l} 
\displaystyle - J_0 \, \frac{h}{\sinh \, (\mathcal{K}\,h) - \mathcal{K}\,h} \,
 \sinh \, \big(\mathcal{K}\, y \big)  + 
J_0   \, \frac{\sinh \, (\mathcal{K}\,h)}{\sinh \, (\mathcal{K}\,h) - \mathcal{K}\,h} 
\, y \ \cosh \big(\mathcal{K}\, y \big) \\  \vspace{-4 mm}  \\ \displaystyle
+ J_0  \, \frac{1-\cosh \, (\mathcal{K}\,h)}{\sinh \, (\mathcal{K}\,h) -\mathcal{K}\,h}
 \, y \, \sinh \big( \mathcal{K} \, y \big)  \,, \qquad  \mathcal{K} = \frac{2 k \pi }{L}  \,. \, 
\end{array} \right. \monend  
The stream function 
$ \, \psi = f(y) \, \cos \big(  \mathcal{K} \, x \big) $, 
the two components $ \, J_x = \frac{\partial \psi}{\partial y} \, $
and  $ \, J_y = - \frac{\partial \psi}{\partial x} \, $
of momentum and the pressure field 
$ \,\, p(x,y)= {{\nu}\over{\mathcal{K}}} \, 
  \sin \,(\mathcal{K}\, x)   
 \big(\frac{{\rm{d}}^3 f}{ {\rm{d}} y^3} - \mathcal{K}^2 \, \frac{{\rm{d}} f}{{\rm{d}} y} \big) \,\, $
define a particular solution of the Stokes problem (\ref{stokes}). 

\monitem 
We have measured the error between the measured values $ \, j_x \, $ 
and $ \, j_y \, $ in the first cell and the four following quantities~:
(i) the given value $\, J_x \, $ and  $\, J_y \, $ on the boundary, 
(ii) the result of the Taylor expansion (\ref{dev_jxjy_o2}) taking into account only 
the first order terms, 
(iii) Taylor expansion (\ref{dev_jxjy_o2}) with all terms of second order
and (iv) the exact values  $\, J_x (x, \, {{\Delta x}\over{2}}) \, $
and  $\, J_y  (x ,\, {{\Delta x}\over{2}}) \, $ of the problem (\ref{stokes})
at the mesh point location. 

\begin{figure}[htbp!]
\begin{center}
\centerline   {\includegraphics[height=.66 \textwidth, width = .88 \textwidth] 
{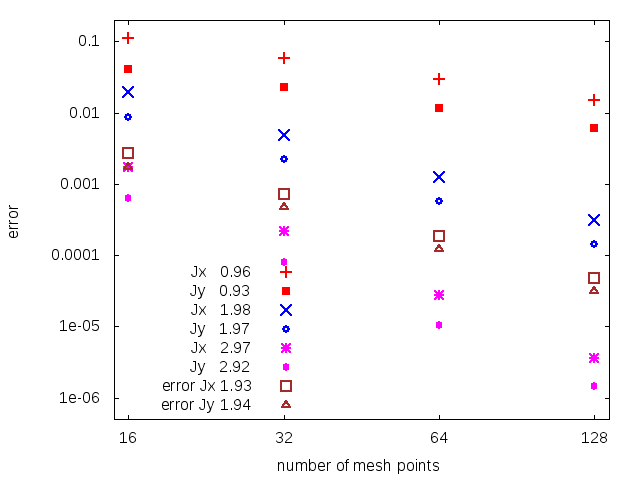}} 
\caption {Theoretical and measured rates of convergence $\theta$  for classical bounce back
(\ref{bb-classique}) as a function of the number of mesh points in the axial direction. 
 The error is proportional to $ \, \Delta x^{\theta} \, $ in the $\ell^{2}$ 
      norm for $x$ and $y$ component of the momentum, for aspect ratio $ \, {{L}\over{h}} = 2 $.}
\label{accordeon-classique} \end{center} \end{figure}

\begin{figure}[H]   
\begin{center}
\centerline   {\includegraphics[height=.66 \textwidth, width = .88 \textwidth] 
{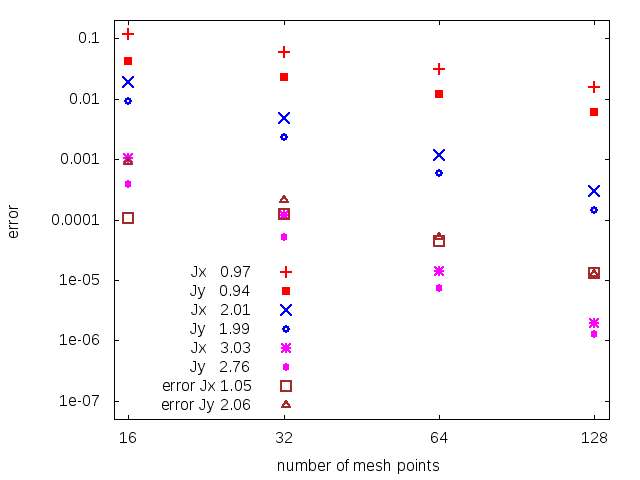}} 
\caption {Theoretical and measured rates of convergence $\theta$  for generalized bounce back
(\ref{BB-O2}) as a function of the number of mesh points in the axial direction. 
 The error is proportional to $ \, \Delta x^{\theta} \, $ in the $\ell^{2}$ 
      norm for $x$ and $y$ component of the momentum, for aspect ratio $ \, {{L}\over{h}} = 2 $.}
\label{accordeon-general} \end{center} \end{figure}

We have done two numerical experiments. 
One (see Fig.~\ref{accordeon-classique})  with very simple values of the coefficients of (\ref{BB-O2})~:
$ \, a_2=a_5=a_6=0, \,  k_2=k_5=k_6=0 \,  $ and the other (see Fig.~\ref{accordeon-general})
with the condition (\ref{choix-parametres-ordre-1}) and the choice 
$ \, a_2=a_5=-1 $, $ \,  k_2=4 $, $ \, k_5=k_6=1  $. 

The results are as  expected. This validates the formal expansion proposed in \cite{DLT15}  for the analysis
of the bounce back boundary condition. 
The error is only first order for  the $x $ component of the momentum. This is due to a particulary 
good precision with only 16 mesh points (see Fig.~\ref{accordeon-general}).

\bigskip {\bf 7) \quad Conclusion}  

We have shown that the classical bounce back is the result of an approximation 
at  order zero of the internal  lattice Boltzmann  scheme.
An analysis by an extension of the  Taylor expansion method  was described as in \cite{DLT15}. 
Then a new scheme called first order bounce back was proposed and analyzed. 
We proved that in this scheme
 the artefact/defect at order 1 of the classical bounce back can be removed.
Finally we proposed an extended bounce back  scheme where we removed all the artefacts/defects 
at order~1 and we proved that
for a special choice the bounce back can be exact up to order two for Poiseuille flow test case. 
The stationary ``accordion'' test case shows that for  a nontrivial flow, 
the analysis proposed for the boundary condition does not present any contradiction.
Other numerical experiments will be presented in forthcoming contributions. 
Moreover, an analysis of the anti-bounce back \cite{Gi05}, appropriate for taking into 
consideration a pressure boundary condition, seems also possible. 

\bigskip {\bf Acknowledgments }  
A part of this work has been realized during the stay of two of us at the 
Beijing  Computational Science Research Center. We thank the colleagues of CSRC 
for their hospitality. Last but not least, the authors  thank the referees  who suggested several points 
 in need of improvement.

\bigskip {\bf References }  
\section*{References}

 \bibliography{mybibfile}

\begin{thebibliography}{99}

 \vspace{-.34cm}  

 \vspace{-.24cm}  
\bibitem{BFL01}  
M. Bouzidi, M. Firdaous, P. Lallemand,  
``Momentum transfer of a Boltzmann-lattice fluid with boundaries'',
\textit {Physics of Fluids}, vol.~{\bf{13}}, p.~3452-3459, 2001.


\vspace{-.24cm}  
\bibitem{DLT10}           
F. Dubois, P. Lallemand, M. M. Tekitek, ``On a superconvergent 
lattice Boltzmann boundary scheme'', 
{\it{Computers and Mathematics with Applications}}, 
vol. {\bf{59}}, p. 2141-2149, 2010.

 \vspace{-.24cm}  
\bibitem{DLT15}   
 F. Dubois, P.  Lallemand, M. M. Tekitek, 
``Taylor expansion method for analysing bounce-back boundary conditions 
for lattice Boltzmann method'', {\it{ESAIM: Proceedings}},  vol.~{\bf 52}, p.~25--46, 2015.      

 \vspace{-.24cm}  
\bibitem{Gi05}   
I. Ginzburg, 
``Generic boundary conditions for lattice Boltzmann models and their appli-
cation to advection and anisotropic dispersion equations'', 
{\it Advances in Water Resources}, vol.~{\bf 28}, p.~1196-1216, 2005.

 \vspace{-.24cm}  
\bibitem{GA94}  
I. Ginzburg, P. Adler, 
``Boundary flow condition analysis for the three-dimensional lattice Boltzmann model'', 
\textit{Journal of Physics II France}, vol.~{\bf 4},  p.~191-214, 1994. 



 \vspace{-.24cm}  
\bibitem{dHG03} 
 D. d'Humi\`eres, I. Ginzburg,  
``Multi-reflection boundary  conditions for lattice Boltzmann models'', 
{\it{Physical Review E}}, vol.~{\bf{68}}, issue~6, p.~066614 (30 pages), 2003.



 \vspace{-.24cm}  
\bibitem{LL00}  
P. Lallemand, L-S. Luo. 
  ``Theory of the lattice Boltzmann method: 
   Dispersion, dissipation, isotropy, Galilean invariance, and stability'',
{\it Physical Review E}, vol.~{\bf 61}, p.~6546-6562, June 2000.  



 \vspace{-.24cm}  
\bibitem{ZH97}  
Q. Zou, X. He. 
``On pressure and velocity boundary conditions for the lattice Boltzmann BGK model'',
{\it {Physics of Fluids}}, vol.~{\bf{9}}, p.~1591-1598, 1997. 



\end{thebibliography}

\end{document}